\documentclass[a4paper,11pt]{article}
\usepackage{amsfonts,amssymb,amsbsy,amsmath,amsthm,mathrsfs,enumerate,verbatim}

  \usepackage{graphics} 
  \usepackage{epsfig} 
\usepackage{graphicx}  \usepackage{epstopdf}
 \usepackage[colorlinks=true]{hyperref}

\hypersetup{urlcolor=blue, citecolor=red}

\usepackage{amsfonts}
\usepackage{amsmath}
\usepackage{amssymb}
\usepackage{graphicx}
\usepackage{epsf}
\usepackage{float}

\usepackage{epsfig}

\usepackage{color}
\usepackage{afterpage}
\usepackage{verbatim}
\usepackage{multicol}
\usepackage{multirow}

\usepackage{times}

\newtheorem{theo}{Theorem}

\theoremstyle{definition}
\newtheorem{definition}{Definition}

\newcommand{\bel}{\begin{equation} \label}
\newcommand{\ee}{\end{equation}}

\def\beq{\begin{equation}}
\def\eeq{\end{equation}}
\newcommand{\bea}{\begin{eqnarray}}
\newcommand{\eea}{\end{eqnarray}}
\newcommand{\beas}{\begin{eqnarray*}}
\newcommand{\eeas}{\end{eqnarray*}}

 \usepackage{graphicx,color}
 \definecolor{mygreen}{cmyk}{1,0,1,0.1}

\usepackage{epsf}
\usepackage{epstopdf}

\usepackage{pdfsync}

\graphicspath{
  {FIGURES/}
}

\begin{document}

\title{Interactive Change Point Detection using
optimisation approach and Bayesian statistics applied to real world
applications}

\author{R. Gedda  \thanks{
Chalmers University of Technology and
ABB Corporate Research Centre, e-mail: \texttt{\
rebeccagedda@outlook.com}}
\and
L. Beilina  \thanks{
Department of Mathematical Sciences, Chalmers University of Technology and
 University of Gothenburg, SE-42196 Gothenburg, Sweden, e-mail: \texttt{\
larisa@chalmers.se}}
\and
R. Tan \thanks{ABB Corporate Research Centre, Ladenburg, Germany, e-mail: \texttt{\
ruomu.tan@de.abb.com}}
}

\date{}
\maketitle
\begin{abstract}
\noindent
\textit{Change point detection} becomes more and more important as datasets increase in size, where unsupervised detection algorithms can help users process data. 
To detect change points, a number of unsupervised algorithms have been developed which are based on different principles. 
One approach is to define an optimisation problem and minimise a cost function along with a penalty function. 
In the optimisation approach, the choice of the cost function affects the predictions made by the algorithm. 
In extension to the existing studies, a new type of cost function using Tikhonov regularisation is introduced.
Another approach uses Bayesian statistics to calculate the posterior probability distribution of a specific point being a change point. 
It uses a priori knowledge on the distance between consecutive change points and a likelihood function with information about the segments.
The optimisation and Bayesian approaches for offline change point detection are studied and applied to simulated datasets as well as a real world multi-phase dataset. 
The approaches have previously been studied separately and a novelty lies in comparing the predictions made by the two approaches in a specific setting, consisting of simulated datasets and a real world example. 
The study has found that the performance of the change point detection algorithms are affected by the features in the data.
%
%
%
\end{abstract}

\maketitle


\section{Motivation}
\label{sec:motivation}
The topic of \textit{Change Point Detection} (CPD) has become more and
more relevant as time series datasets increase in size and often
contain repeated patterns. To detect change points in data,
segmentation can be performed to group similar phases of the time
series data together. This is of importance for complicated and large
datasets where exist multiple phases which are desirable to separate
in order to compare data.  To detect the change points, a number of
algorithms have been developed which are based on different
principles \cite{continuous_inspection_schemes, E.S.Page_2,
selective_review_CPD, wireless_sensros, financial_CPD,
traffic_model_cpd, CPD_chemical}.
\\
\\
The first work on change point detection was done by
Page \cite{continuous_inspection_schemes, E.S.Page_2} where piecewise
identically distributed datasets were studied. The objective was to
identify various features in the independent and non-overlapping
segments. Examples of features can be mean, variance and distribution
function for each data segment.  Detection of change points can either
be done in real time or in retrospective, and for a single signal or
in multiple dimensions. The real time approach is generally known as
online detection, while the retrospective approach is known as offline
detection. This work is based on offline detection, meaning all data
is available for the entire time interval under investigation. Many
CPD algorithms are generalised for usage on multi-dimensional
data \cite{selective_review_CPD} where one-dimensional data can be
seen as a special case. This work focuses on one-dimensional time
dependent data, where results are more intuitive and common in real
world settings. Another important assumption is connected to the
number of change points in the data. This can either be known
beforehand or unknown. This work assumes that the number of change
points is not known. Various CPD methods have been applied to a vast
spread of areas, stretching from sensor
signals \cite{wireless_sensros} to natural language
processing \cite{language_processing}.  Some CPD methods have also
been implemented for financial analysis \cite{financial_CPD} and
network systems \cite{traffic_model_cpd}, where the algorithms are
able to detect changes in the underlying setting.  Change point
detection has also been applied to chemical
processes \cite{CPD_chemical}, where the change points in the mean of
the data from chemical processes are considered to be the
representation of changed quality of production.  This illustrates the
usability for change point detection and presents the need of domain
expert knowledge to confirm that the algorithms make correct
predictions.
\\
\\
The current work is based on numerical testing of two approaches for
CPD: the optimisation approach, with and without regularisation, and the
Bayesian approach.
Both approaches are tested on real world data from a multi-phase
flow facility. These  approaches are developed and studied separately in
previous studies \cite{selective_review_CPD,
Exact_efficient_Bayesian_inference}. Performance comparison and the
evaluation of computational efficiency of both approaches form the
foundation of this work. The work by Truong et
al. \cite{selective_review_CPD} gives a good overview of change point
detection algorithms which are based on the optimisation approach.  In
extension to the existing work, new cost functions based on
regularisation can be implemented. Some examples where regularisation
techniques are used in machine learning algorithms are presented
in \cite{pattern_recognition_and_ML, supervised_ML:_Classification,
Goodfellow}.  For the Bayesian approach, the work by
Fearnhead \cite{Exact_efficient_Bayesian_inference} gives a thorough
description of the mathematics behind the algorithm. These two
approaches have been studied separately and this work aims to compare
the predictions made by the two approaches in a specific setting of
simulated data and a real world example. For numerical comparison, the methods
presented by van den Burg and Williams \cite{Turing_paper} are used,
along with metrics specified by Truong et
al. \cite{selective_review_CPD}. Throughout the work, the following
two main questions are considered:
\\
\begin{itemize}

\item \textit{How are the two change point detection  approaches 
affected by the features of investigated data?}

\item \textit{How can user knowledge and feedback be incorporated in the two above mentioned approaches?}

\end{itemize}
\vspace{0.2cm}
The first question formulates the main research
investigation of this work and suggests that studied algorithms should be compared. The secondary focus of the work lies in the domain expert interplay, where the possibility of interaction is studied.
\\
\\
This work is structured as follows: first the appropriate notation is introduced along with definitions. Then, the two approaches, optimisation and Bayesian, are derived separately in section \ref{sec:background} along with the metrics used for comparison. Section \ref{sec:methods} describes the used datasets and the testing procedure. The results of the study are presented in section \ref{sec:results} and a discussion is held in section \ref{sec:discussion}. Finally, a summary of findings and conclusions is provided in section \ref{sec:summary}.

\section{Background}
\label{sec:background}
In this section we provide the necessary background knowledge for the exploration and compare the two CPD approaches. First, the notation used in the entirety of this work is introduced and definitions for change points are provided. The optimisation approach's components, penalty function, search direction and cost functions, are explored. The Bayesian approach is derived from Bayes' formula to the problem specific formulation used in this work. Finally, the test metrics used for evaluation are introduced. 

\subsection{Notation and definition}
\begin{figure}
    \centering
    \includegraphics[width = 1 \linewidth]{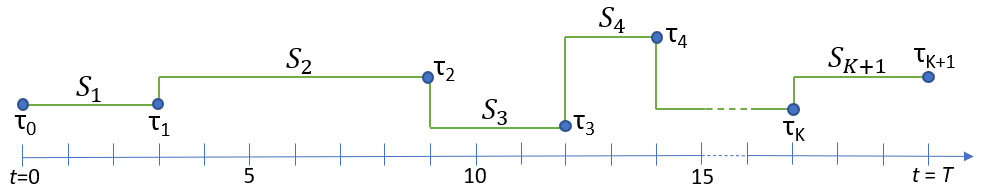}
    \caption{Illustration of used notation. In the figure, we see how $K$ intermediate change points are present on a time interval $t \in [0,T]$, where we note how $\tau_0$ and $\tau_{K+1}$ are synthetic change points. }
    \label{fig:break_point_notation}
\end{figure}

Figure~\ref{fig:break_point_notation} shows multiple change points ($\tau_0, \tau_1, ..., \tau_{K+1}$) and the segments ($S_1, S_2, ..., S_{K+1}$), defined by the change points, and is an example of a uni-variate time-series dataset. The purpose of CPD is to find the time stamps $\tau_0$ to $\tau_{K+1}$ given the time series.
Let us now introduce the notations used in the work.
\\
\\
Throughout the work we are working in the time domain $[0,T] \subset \mathbb{N}$ which is discretised with $t_i \in [0,T], \; i = 0, 1, ... , n$. The signal value at time $t_i$ is given by $y_{t_i} := y(t_i)$.
The time points are equidistant, meaning $t_{i+1} = t_i + dt$, with $dt = \frac{T}{n-1},\; \forall i$. Let us denote the jump by $[y_{t_i}]$ in time of the discrete function $y_{t_i}$ at time moment $t_i$ which we define as
\begin{equation}\label{jump}
    [y_{t_i}] :=  \lim_{s \to 0+} (y(t_{i+s}) - y(t_{i-s})).
\end{equation}
A set of $K$ change points is denoted by $\mathcal{T}$ and is a subset of time indices $\{ 1,2, \dots,  n\}$. 
The individual change points are indicated as $\tau_j$, with $j \in \{0,1,\dots, K, K+1 \}$, where $\tau_0 = 0$ and $\tau_{K+1} = T$.
Note, with this definition the first and final time points are implicit change points, and we have $K$ intermediate change points ($|\mathcal{T}| = |\{ \tau_1, \dots , \tau_K \}| = K$). 
A segment of the signal from $a$ to $b$ is denoted as $y_{a : b}$, where $y_{0 : T}$  means  the entire signal. With the introduced notation for change points, the segment $S_j$ between change points $\tau_{j-1}$ and $\tau_j$ is defined as  $S_j := [\tau_{j-1},  \tau_{j}], |S_j| = \tau_j - \tau_{j-1}$.
$S_j$ is the $j$-th non-overlapping segment in the signal, $j \in \{1, \dots, K+1 \}$, see Figure~\ref{fig:break_point_notation}. Note that the definition for $S_j$ does not hold for $j=0$, since this is the first change point. 
\\
\\
As the goal is to identify whether a change has occurred in the signal, a proper definition of the term change point is needed along with clarification of change point detection. Change point detection is closely related to change point estimation (also known as change point mining, see \cite{change_point_mining, change_analysis}).
According to Aminikhanghahi and Cook~\cite{survey_of_methods_in_CPD}, change point estimation tries to model and interpret known changes in time series, while change point detection tries to identify whether a change has occurred \cite{survey_of_methods_in_CPD}. This illustrates that we will not focus on the change points' characteristics, but rather if a change point exists or not.
\\
\\
One challenge is to identify the number of change points in a given time series. The problem is a balance between having enough change points whilst not over-fitting to the data. If the number of change points is known beforehand the problem is merely a best fit problem. 
On the other hand, if the number of change points is not known, the problem can be seen as an optimisation problem with a penalty term for every added change point, or as enforcing a threshold when we are certain enough that a change point exists.
It is evident that we need clear definitions of change points in order to detect them. 
The definitions for change point and change point detection are defined below and will be used throughout this work.

\begin{definition}[Change point]
A change point represents a transition between different states in a signal or dataset.
If two consecutive segments $y_{t_l:t_i}$ and $y_{t_{i},t_m}$, $t_l< t_i < t_m,\; \; l, i, m = 0,...,n$, defined as
  \begin{equation}\label{cpl}
    \begin{split}
  y_{t_l:t_i}  &= \{ t_k | [y(t_k)] < |\varepsilon_l|,\; l \leq k \leq i \}, \\
  y_{t_i:t_m}  &= \{ t_k | [y(t_k)] < |\varepsilon_m|,\; i\leq  k \leq m \},
  \end{split}
  \end{equation}
have a distinct change in features such that
 $|\varepsilon_l| << |\varepsilon_m|$ or $|\varepsilon_l| >> |\varepsilon_m|$, or if $y_{t_{i}}$ is a local extreme point (i.e minimum or maximum\footnote{If $f(x^*) \leq f(x)$ or $f(x^*) \geq f(x)$ for all $x$ in $X$ within distance $\epsilon$ of $x^*$, then $x^*$ is a local extreme point.})
then, $\tau_j = t_i,\; j= 0,1,...,K,K+1$ is a change point between the two segments.
\label{def:CP}
\end{definition}
Remark: the first change point $\tau_0 = t_0$ and the final change point $\tau_{K+1} = t_n$ are artificial change points which are used to define segments. These two points are defined the same for all predictions and are not part of the prediction process. 
We note that the meaning of a \emph{distinct change} in this definition is different for different CPD methods, and it is discussed in detail in section 3. 
This definition is useful when dealing with features in the data, but there may be other types of change points. These points can be more complex to identify, but are of interest for a domain expert. These change points are referred to as domain specific change points and are defined below.

\begin{definition}[Change point, domain specific] \label{def:CP_domain}
For some process data, a change point is where a phase in the process starts or ends. These points can be indicated in the data, or be the points in the process of specific interest without a general change in data features. 
\end{definition}
\noindent
Finally, we give one more definition of CPD for the case of available information about the probability distribution of a stochastic process.

\begin{definition}[Change point detection]
Identification of times when the probability distribution of a stochastic process or time series segment changes. This concerns detecting whether or not a change has occurred, or whether several changes might have occurred, and identifying the times of any such changes.
\label{def:CPD}
\end{definition}

\subsection{Optimisation approach}
\begin{figure}
    \centering
    \includegraphics[width = 0.8\linewidth]{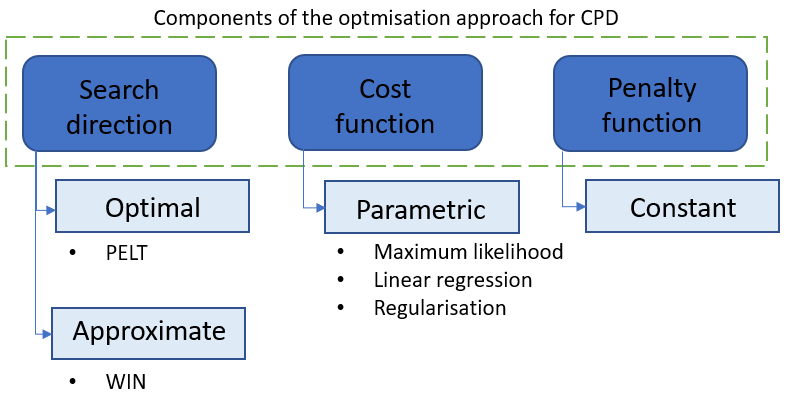}
    \caption{Illustration of the components used in the optimisation approach for change point detection. Each of the three components illustrate the strategies studied in this work. }
    \label{fig:CPD_components}
\end{figure}

Solving the task of identifying change points in a time series can be done by formulating an optimisation problem. A detailed presentation of the framework is given in the work by Truonga et al \cite{selective_review_CPD}, while only a brief description is presented here. The purpose is to identify all the change points, without detecting fallacious ones. Therefore, the problem is formulated as a minimisation problem, where we strive to minimise the cost of segments and penalty per added change point. We need this penalty since we do not know how many change points will be presented. Mathematically, the non-regularised optimisation problem is formulated as
\begin{equation}
    \min_\mathcal{T} \; V(\mathcal{T}) \; + \; pen(\mathcal{T}) \; = \;
    \min_\mathcal{T} \; \sum_{j = 1}^K c(y_{\tau_j:\tau_{j+1}}) \; + \; \beta |\mathcal{T}|,
    \label{eq:opt_problem}
\end{equation}
while the regularised analogy is
\begin{equation}
    \min_\mathcal{T} \; V(\mathcal{T}) \; + \; pen(\mathcal{T})  + reg(\mathcal{T}) \; = \;
    \min_\mathcal{T} \; \sum_{j = 1}^K c(y_{\tau_j:\tau_{j+1}}) \; + \; \beta |\mathcal{T}| + \gamma [\mathcal{T}].
    \label{eq:opt_problem_reg}
\end{equation}
\noindent
Here, $V(\mathcal{T})$ represents a cost function, $pen(\mathcal{T}) $ is a linear penalty function with constant $\beta$ and $reg(\mathcal{T})$ is a regularisation term in appropriate norm in the time space $[0,T]$ with the regularisation parameter $\gamma$.
\\
\\
To solve the optimisation problem \eqref{eq:opt_problem}, we need three components to be combined together: the search method, the cost function and the penalty term.
Figure~\ref{fig:CPD_components} shows a schematic view of how search method, cost function and penalty term create components of a CPD algorithm. 
There are numerous combinations of components which can be chosen for problem \eqref{eq:opt_problem}. Figure~\ref{fig:CPD_components} also illustrates which methods will be studied in this work. The two methods for search directions and cost functions are presented in separate sections, while the choice of penalty function is kept brief. A common choice of penalty function is a linear penalty, which means each added change point $\tau_j$ corresponds to a penalty of $\beta$. A summary of other combinations are presented in Table 2 in the work by Truonga et al. \cite{selective_review_CPD}.

\subsubsection{Search direction}
\label{sec:opt_search_dir}
The search method poses a trade-off between accuracy and computational complexity. In CPD
there are two main approaches used for this, optimal and approximate, see Figure~\ref{fig:CPD_components}.
\\
\\
The problem formulated in equation \eqref{eq:opt_problem} should be solved for an unknown $K$, where the penalty function can be chosen as a constant function, $pen(\cdot) = \beta$. The search method used for this special case is known as \textit{Pruned Exact Linear Time} (abbreviated \texttt{PELT}) and implements a pruning rule. The pruning rule states that for two indices $s$ and $t$, $s<t<T$, if the following condition holds
\begin{equation*}
    \Big[ \min_{\mathcal{T}} \; V(\mathcal{T}, y_{0:s}) + \beta |\mathcal{T}| \Big] \; + \; c(y_{s:t}) \geq \Big[ \min_{\mathcal{T}} \; V(\mathcal{T}, y_{0:t}) + \beta |\mathcal{T}| \Big],
\end{equation*}
then $s$ cannot be the last change point. Intuitively, the algorithm compares if it is beneficial to add another change point between $s$ and $t$. If the cost of a segment $y_{s:t}$ is greater than the cost of two separated segments $y_{s:\tau}$, $y_{\tau+1:t}$ and the additional penalty $\beta$, then there is a change point $\tau$ present between indices $s$ and $t$. The \texttt{PELT}-algorithm is presented in Algorithm 1 in \cite{masterThesis}, and has a time complexity $\mathcal{O}(T)$ \cite{selective_review_CPD}. A drawback of this algorithm is that it can become computationally expensive for large datasets with many time stamps $t$. 
\\
\\
An alternative approach is to use an approximate search direction algorithm to reduce complexity. To reduce the number of performed calculations, an approximate search direction can be used, where partial detection is common. A frequently used technique is the \textit{Window-sliding} algorithm (denoted as \texttt{WIN}-algorithm), when the algorithm returns an estimated change point in each iteration. Similar to the concept used in the \texttt{PELT}-algorithm, the value of the cost function between segments are compared. This is known as the discrepancy between segments and is defined as
\begin{equation*}
    Disc(y_{t-w:t}, y_{t:t+w}) \; = \; c(y_{t-w:t+w}) - \big( c(y_{t-w:t}) + c(y_{t:t+w}) \big),
\end{equation*}
where $w$ is defined as half of the window width. Intuitively, this is merely the reduced cost of adding a change point at $t$ in the middle of the window. The discrepancy is calculated for all $w \leq t \leq T-w$. When all calculations are done, the peaks of the discrepancy values are selected as the most profitable change points. The algorithm is provided in Algorithm 2 in \cite{masterThesis}. There are other approximate search directions, which are not covered in this work, presented by Trounga et al. \cite{selective_review_CPD}.
For this work, the \texttt{PELT}-algorithm is used for the optimal approach and the \texttt{WIN}-algorithm is used for the approximate approach.



\subsubsection{Cost functions}
\label{sec:opt_cost_functions}
The cost function can decide which feature changes are detected in the data. In other words, the cost function measures the homogeneity. There are two approaches for defining a cost function; parametric and non-parametric. The respective approaches assume either that there is an underlying distribution in the data, or that there is no distribution in the data. This work focuses on the parametric cost functions, for three sub-techniques illustrated in Figure~\ref{fig:CPD_components}. The three techniques, maximum likelihood estimation, linear regression and regularisation, are introduced in later sections with corresponding cost function definitions. 
\\
\\
\textit{Maximum Likelihood Estimation} (MLE) is a powerful tool with a wide application area in statics.
MLE finds the values of the model parameters that maximise the likelihood function $f(y|\Theta)$ over the parameter space $\mathbf{\Theta}$ such that
\begin{equation*}
MLE(y) = \max_{\Theta \in  \mathbf{\Theta}  } f(y|\Theta),
\end{equation*}
where $y$ is observed data and $\Theta \in  \mathbf{\Theta}$ is a vector of parameters. In the setting of change point detection, we assume the samples are independent random variables, linked to the distribution of a segment. This means that for all $t\in[0,T]$, the sample
\begin{equation}
    y_t \sim \sum_{j = 0}^{K} f(\cdot | \theta_j) \; \mathbf{1}(\tau_j < t < \tau_{j+1}),
    \label{eq:y_t_distribution}
\end{equation}
where $\theta_j$ is a segment specific parameter for the distribution. The function $\mathbf{1(\cdot)}$ is the delta function $\delta([\tau_j,\tau_{j+1}])$, and is equal to one if sample $y_t$ belongs to segment $j$, otherwise zero:
\begin{equation*}
  \mathbf{1}(\tau_j <  t < \tau_{j+1}):= \delta([\tau_j,\tau_{j+1}]) = \left\{ \begin{array}{ll} 1 \quad  & \text{if }
     y_t \in [\tau_j,\tau_{j+1}], \\
   0\qquad & \text{elsewhere.} \end{array} \right. 
\end{equation*}
The function $f(\cdot | \theta_j) $  in \eqref{eq:y_t_distribution}
represents the likelihood function for the distribution with parameter $\theta_j$.
Then the $ MLE(y_t)$    reads:
\begin{equation*}
MLE(y_t) = \max_{\theta_j \in  \mathbf{\Theta} } \sum_{j = 0}^{K} f(\cdot | \theta_j) \; \mathbf{1}(\tau_j < t < \tau_{j+1}) 
\end{equation*}
where $\theta_j$ is segment specific parameter for the distribution. Using $MLE(y_t)$ we can estimate the segment parameters $\theta_j$, which are the features in the data that change at the change points. If the distribution family of $f$ is known and
the sum of costs, $V$ in \eqref{eq:opt_problem} or \eqref{eq:opt_problem_reg}, is equal to the negative log-likelihood of $f$, then MLE is equivalent to change point detection. Generally, the distribution $f$ is not known, and therefore the cost function cannot be defined as the negative log-likelihood of $f$. 
\\
\\
In some datasets, we can assume the segments to follow a Gaussian distribution, with parameters mean and variance.
More precisely, if $f$ is a Gaussian distribution, the MLE for expected value (which is the distribution mean) is the sample mean.
If we want to identify a shift in the mean between segments, but where the variance is constant, the cost function can be defined as the quadratic error between a sample and the MLE of the mean. For a sample $y_t$ and the segment mean $\Bar{y}_{a:b}$ the cost function is defined as
\begin{equation}
    c_{L2}(y_{a:b}) \; := \; \sum_{t = a+1}^b || y_t - \Bar{y}_{a:b} ||_2 ^2,
    \label{eq:c_l2}
\end{equation}
where  the norm  $\| \cdot \|_2$ is the usual $L_2$-norm defined for any vector $v \in \mathbb{R}^n$  as
$$
\|v\|_{2}:={\sqrt {(v_{1})^{2}+(v_{2})^{2}+\dotsb +(v_{n})^{2}}}.
$$
The cost function \eqref{eq:c_l2} can be simplified for uni-variate signals to 
$$ c_{L2}(y_{a:b}) := \sum_{t = a+1}^b ( y_t - \Bar{y}_{a:b})^2 $$
which is equal to the MLE variance times length of the segment. More explicitly, for the presumed Gaussian distribution $f$ the MLE of the segment variance $\Hat{\sigma}^2_{a:b}$ is calculated as $\Hat{\sigma}^2_{a:b} = \frac{c_{L2}(y_{a:b})}{b-a}$, using the MLE of the segment mean, $\Bar{y}_{a:b}$. This estimated variance $\Hat{\sigma}^2_{a:b}$ times the number of samples in the segment is used as the cost function for a segment $y_{a:b}$. This cost function is appropriate for piecewise constant signals, shown in Figure~\ref{fig:break_point_notation}, where the sample mean $\Bar{y}_{a:b}$ is the main parameter which changes. We note that this formulation mainly focuses on changes in the mean, and the cost is given by the magnitude of the variance of the segment around this mean. 
A similar formulation can be given in the $L_1$-norm,
\begin{equation}
    c_{L1}(y_{a:b}) \; := \; \sum_{t = a+1}^b | y_t - \Tilde{y}_{a:b} |,
    \label{eq:c_l1}
\end{equation}
where we find the least absolute deviation from the median $\Tilde{y}_{a:b}$ of the segment. Similar to the cost function in equation \eqref{eq:c_l2}, the cost is calculated as the aggregated deviation from the median for all samples in $y_{a:b}$. This uses the MLE of the deviation in the segment, compared to the MLE estimation of the variance used in \eqref{eq:c_l2}. Again, the function mainly identifies changes in the median, as long as the absolute deviation is smaller than the change in median between segments. 
\\
\\
An extension of cost function \eqref{eq:c_l2} can be made to account for changes in the variance.
The empirical covariance matrix $\hat{\Sigma}$ can be calculated for a segment from $a$ to $b$. The cost functions for multi- and uni-variate signals are defined
by \eqref{eq:c_normal_mulivar} and  \eqref{eq:c_normal_univar},
correspondingly, as
\begin{equation}
    c_{Normal}(y_{a:b}) \; := \; (b-a) \log \: \det \: \hat{\Sigma}_{a:b} \; + \; \sum_{t = a+1}^b ( y_t - \Bar{y}_{a:b} )'\: \hat{\Sigma}_{a:b}^{-1} \:( y_t - \Bar{y}_{a:b} ),
    \label{eq:c_normal_mulivar}
\end{equation}
\begin{equation}
   c_{Normal}(y_{a:b}) \; :=  \; (b-a) \log \: \hat{\sigma}^2_{a:b} \; + \; \frac{1}{\hat{\sigma}^2_{a:b}}\sum_{t = a+1}^b \; ( y_t - \Bar{y}_{a:b} )^2,
    \label{eq:c_normal_univar}
\end{equation}
where $\hat{\sigma}_{a:b}$ is the empirical variance of segment $y_{a:b}$. For the uni-variate case, we note that
$$
c_{\Sigma}(y_{a:b}) = (b-a) \log \: \hat{\sigma}^2_{a:b} +\hat{\sigma}^{-2}_{a:b}
c_{L_2}(y_{a:b}),
$$
which clearly is an extension of equation \eqref{eq:c_l2}. This cost function is appropriate for segments that follow Gaussian distributions, where both the mean and variance parameters change between segments. 
\\
\\
If segments in the signal follow a linear trend, a linear regression model can be fitted to the different segments. At change points, the linear trends in the respective segment changes abruptly. In contrast to the assumption formulated in \eqref{eq:y_t_distribution}, the assumption for linear regression models is formulated as
\begin{equation*}
    y_{t}=\mathbf {x} _{t}^{\mathsf {T}}{\boldsymbol {\beta }}+\varepsilon _{t} \; = \; \beta_j^0 \: + \: \beta_j \: x_t \; + \; \varepsilon_t, \quad \forall t, \: \tau_j < t < \tau_{j+1}
\end{equation*}
with the intercept $\beta_j^0$ and coefficient $\beta_j$ dependent on segment $j=\{ 0, \dots , K+1\}$. The noise for each sample is given by $\varepsilon_t$, which is assumed to be normally distributed with zero mean. Having only one covariate $x_t$, the model fitting is known as a simple linear regression model, which constitutes an intercept and a coefficient for the covariate. 
The intercept $\beta_j^0$ and coefficient $\beta_j$ are unknown and each segment is presumed to have an underlying linear regression model. A simple minimisation problem for the cost function which uses the simple linear regression is defined as
\begin{equation}
    c_{LinReg}(y_{a:b}) \; := \; \min_{\beta \in \mathbb{R}^p} \; \sum_{t = a+1}^b \Big( y_t \: - \:(\beta_j^0 + \:\beta x_t)  \Big)^2,
    \label{eq:c_linear}
\end{equation}
where we use a single covariate $x_t$.
The cost is given by the error between the simple linear regression and the samples, and is known as the model squared residual.
\\
\\
If we use previous samples $[y_{t-1}, y_{t-2},..., y_{t-p}]$ as covariates, we have an autoregressive model. In this thesis, this is limited to four lags ($p=4$), meaning the covariate at $t = t_i$ is defined as the vector $\mathbf{\Tilde{x}}_{t_i} = [y_{{t_i}-1}, y_{{t_i}-2}, y_{{t_i}-3}, y_{{t_i}-4}]$. Similar to equation~\eqref{eq:c_linear}, we can define a cost function as
\begin{equation}
    c_{AR}(y_{a:b}) \; := \; \min_{\beta \in \mathbb{R}^p} \; \sum_{t = a+1}^b || y_t \: - \: (\beta_0 + \:\beta \mathbf{\Tilde{x}}_{t})  ||_2^2,
    \label{eq:c_AR}
\end{equation}
where $\mathbf{\Tilde{x}}_{t}$ is a collection of $p$ lagged samples of $y_t$. This formulation can detect changes in models applied to non-stationary processes. 
\\
\\
By adding a regularisation term to equation \eqref{eq:c_linear}, we can add information and penalise over-fitting. The regularisation term is dependent on the model parameters $\beta$ and a regularisation parameter $\gamma$, where $\gamma$ can be estimated or chosen as a constant ($ \gamma > 0$). If $\gamma = 0$, we get the ordinary linear regression model, presented in equation \eqref{eq:c_linear}. The use of regularisation has been studied widely, where the approach can provide a theoretical, numerical or iterative solution for ill-posed problems \cite{tikhonov_solution_of_ill_posed_problems, tikhonov_numerical_methods_for_ill_posed_problems, iter_regularizartion_for_ill_posed_problems}. 
Tikhonov's regularisation has been used when solving inverse problems \cite{iter_methods_for_approx_solution_Inverse_problems, Approx_global_conv_and_adapt_for_Inverse_problems}
and in machine learning for classification and pattern recognition, see details and analysis of different methods in \cite{pattern_recognition_and_ML, Goodfellow, supervised_ML:_Classification, num_analysis_LS_and_perceptron_for_classification}. In this thesis we study Ridge and Lasso regularisation which are standard approaches of Tikhonov regularisation \cite{ridge_regularisation,lasso_regularisation}.
\\
\\
The first regularisation approach which is studied in this thesis is the Ridge regression,
\begin{equation}
    c_{Ridge} (y_{a:b}) \; := \;  \min_{\beta \in \mathbb{R}^p } \; \sum_{t = a+1}^b \Big( y_t \: - \: (\beta_0 + \:\beta x_t)  \Big)^2 \; + \; \gamma \sum_{j = 1}^p || \beta_j ||_2^2,
    \label{eq:c_ridge}
\end{equation}
where the regularisation term is the aggregated squared $L_2$-norm of the model coefficients. If the $L_2$-norm is exchanged for the $L_1$-norm we get Lasso regularisation. The cost functions is defined as
\begin{equation}
    c_{Lasso} (y_{a:b}) \; := \; \min_{\beta \in \mathbb{R}^p } \; \sum_{t = a+1}^b \Big( y_t \: - \: (\beta_0 + \:\beta x_t)  \Big)^2 \; + \;  \gamma \sum_{j = 1}^p |\beta_j|
    \label{eq:c_lasso}
\end{equation}
where $\gamma$ is the previously described regularisation parameter. Note that this parameter can be the same as the parameter in the Ridge regression \eqref{eq:c_ridge} but these are not necessarily equal. 

\subsection{Bayesian approach}
\label{sec:bayes}
In contrast to the optimisation approach, the Bayesian approach is based on Bayes' probability theorem, where the maximum probabilities are identified.
It is based on the Bayesian principle of calculating a posterior distribution of a time stamp being a change point, given a prior and a likelihood function. From this posterior distribution, we can identify the points which are most likely to be change points. The upcoming section will briefly go through the theory behind the Bayesian approach. For more details and proofs of used Theorem, the reader is directed to the work by Fearnhead  \cite{Exact_efficient_Bayesian_inference}. The section is formulated as a derivation of the sought after posterior distribution for the change points. Using two probabilistic quantities $P$ and $Q$, we can rewrite Bayes' formula to a problem specific formulation which gives the posterior probability of a change point $\tau$. Finally, we combine the individual posterior distribution to get a joint distribution for all possible change points. 
\\
\\
The principle behind the Bayesian approach lies in the probabilistic relationship formulated by Bayes in 1976 \cite{BayesLetter}, where a posterior probability distribution can be expressed as
\begin{equation}
   \mathbf{Pr}(a | b) = \frac{\mathbf{Pr}(a, b)}{\mathbf{Pr}(b)} = \frac{\mathbf{Pr}(b | a) \; \mathbf{Pr}(a)}{\mathbf{Pr}(b)}
   \label{eq:Bayes_formula}
\end{equation}
for the event $a$ given another event $b$. Here, $\mathbf{Pr}(b | a)$ is the likelihood of $b$ given $a$. The distribution $\mathbf{Pr}(a)$ is known as the prior distribution of $a$. 
As soon as $\mathbf{Pr}(b|a)$ and $\mathbf{Pr}(a)$  are defined, the estimator of the posterior distribution $\mathbf{Pr}(a|b)$ can be calculated. A common technique is the \textit{Maximum A Posteriori} (MAP) approach, which is the solution of the problem
\begin{equation*}
MAP(b) = \max_{a \in A} \; \mathbf{Pr}(a|b),
\end{equation*}
where $A$ are the possible values for $a$.
Taking the log of the  above equation, we get
\begin{equation}
  \max_{a \in A} \; \log \mathbf{Pr}(a|b) =\max_{a \in A} [\; \log \mathbf{Pr}(b|a) + \log  \mathbf{Pr}(a) - \log \mathbf{Pr}(b) \;] 
  \label{MAP}
\end{equation}
which is used in this work. 
\\
\\
In our case, we wish to predict the probability of a change point $\tau$ given the data  $y_{1:n}$. Thus, Bayes' formula in \eqref{eq:Bayes_formula} can be reformulated for our problem as
\begin{equation}
    \mathbf{Pr}(\tau | y_{1:n}) = \frac{\mathbf{Pr}( \tau, y_{1:n})}{\mathbf{Pr}(y_{1:n})} = \frac{\mathbf{Pr}(y_{1:\tau} | \tau) \; \mathbf{Pr}(y_{\tau+1:n} | \tau) \; \mathbf{Pr}(\tau)}{\mathbf{Pr}(y_{1:n})},
    \label{eq:symbolic_posterior}
\end{equation}
where $\mathbf{Pr}(y_{1:\tau} | \tau)$ and $\mathbf{Pr}(y_{\tau+1:n} | \tau)$ are the likelihood of segments before and after the given change point $\tau$, respectively. The prior distribution $\mathbf{Pr}(\tau)$ indicates the probability of a potential change point $\tau$ existing and $y_{1:n}$ represents the entirety of the signal. Using the MAP(y) in logarithmic terms, we get the problem specific version of \eqref{MAP}
$$
\max_{\tau \in [0:n]} \log \mathbf{Pr}(\tau | y_{1:n}) =\max_{\tau \in [0:n]} [\log \mathbf{Pr}(y_{1:\tau} | \tau) + \log \mathbf{Pr}(y_{\tau+1:n} | \tau) + \log  \mathbf{Pr}(\tau) - \log \mathbf{Pr}( y_{1:n})].
$$
\\
Similarly to Fearnhead \cite{Exact_efficient_Bayesian_inference}, we will define two functions $P$ and $Q$ which are used for calculations in the Bayesian approach. First, we define the probability of a segment $P(t,s)$, given two entries belonging to the same segment
\begin{equation}
    P(t,s) = \mathbf{Pr}(y_{t:s} | \text{$t,s \in S_j$}) = \int \prod_{i = t}^s \; f(y_i | \theta_{S_j} ) \; \pi(\theta_{S_j}) \; d\theta_{S_j},
    \label{eq:P}
\end{equation}
where $f$ is the probability density function of entry $y_t$ belonging to a segment $S_j$ with parameter $\theta_{S_j}$. We note that this function has similarities used in the optimisation approach, namely in equation \eqref{eq:y_t_distribution}, where we assume a distribution for each segment. In this work, this likelihood will be the Gaussian observation log-likelihood function, but other function choices can be made. Note that $\pi(\theta_{S_j})$ is the prior for the parameters of segment $S_j$. The discrete intervals $[t,s], \; t\leq s$ makes $P$ an upper triangular matrix which elements are probabilities for segments $y_{t:s}$. Note that this probability is independent of the number of true change points $K$.
\\
\\
The second function, $Q$, indicates the probability of a final segment $y_{t:n}$ starting at time $t_i$ given a change point at previous time step, $t_{i-1}$. This probability is affected by the number of change points $K$, and also which of the change points that is located at time $t_{i-1}$. Since we do not know the exact number of change points $K$, we use a generic variable $k$, and perform calculations for all possible values of $K$. The recurrent function is defined as
\begin{align}
    Q_j^{(k)}(i) &=  \mathbf{Pr}(y_{i:n} | \text{$\tau_j = i-1$, $k$}) = \nonumber\\
    &= \sum_{s = i}^{n-k+j} \; P(t,s) \; Q_{j+1}^{(k)}(s+1) \; \pi_k(\tau_j = i-1 | \tau_{j+1} = s),
    \label{eq:Q}
    \\
    Q^{(k)}(1) &= \mathbf{Pr}(y_{1:n} | k) = \nonumber\\
    &= \sum_{s = 1}^{n-k} P(1,s) \; Q_{1}^{(k)}(s+1),
    \label{eq:Q1}
\end{align}
where $Q^{(k)}(1)$ is the first time step and is a special case of $Q_j^{(k)}(i)$. The time index is indicated with $i \in [2, \dots , n]$. The assumed number of change points is denoted $k \in \{1, \dots , n-1\}$, where $j \in \{1, \dots, k\}$ indicates which of the $k$ assumed change points we are currently at. The prior $\pi_k$ is based on the distance between change points, naturally dependent on $k$. This prior can be any point process, where the simplest example is the constant prior with probability $p=1/n$, where $n$ is the number of samples. Other examples include the negative binomial and Poisson distribution. Note that the prior should be a point process since we have discrete time steps. The first time step is defined as an altered function in \eqref{eq:Q1}. The result from this recursion is saved in an array of length $n$. A derivation and proof for this function $Q$ is provided by Fearnhead in Theorem 1 \cite{Exact_efficient_Bayesian_inference}. When calculating the sums in equations \eqref{eq:Q}-\eqref{eq:Q1}, the terms on the right hand side contribute to the function value. We can implement a truncation, with negligible error, at the $k$-th term if 
$$
\frac{P(t,k)\; Q(s+1) \; \pi(k+1-t)}{\sum_{s=t}^{k} P(t,s)\; Q(s+1) \; \pi(s+1-t) } < \epsilon,
$$
where $\pi$ represents the prior distribution for the distance between two consecutive change points and $\epsilon$ is a truncation threshold.
 In this work was used $\epsilon = 10^{-10}$ as a truncation threshold. 
\\
\\
Using functions $P$ and $Q$, the posterior distribution $\mathbf{Pr}(\tau_j | \tau_{j-1}, y_{1:n}, k )$ for change point $\tau_j$, given the previous change point $\tau_{j-1}$, the data $y_{1:n}$ and number of change points, can be calculated.  Using equation~\eqref{eq:symbolic_posterior} along with the expressions for $P$ and $Q$, we can formulate the posterior distribution  for change point $\tau_j$ as
\begin{align}
    \mathbf{Pr}&(\tau_j | \tau_{j-1}, y_{1:n}, k ) = \label{eq:posteriour_tau_j} \nonumber \\
    &= \frac{ \mathbf{Pr}(y_{\tau_{j-1}+1: \tau_j} | \tau_{j-1}+1, \tau_j \in S_j)\; \mathbf{Pr}(y_{\tau_j + 1:n} | \text{$\tau_j$, $k$}) \; \pi_k(\tau_{j-1}| \tau_j) }{\mathbf{Pr}(y_{\tau_{j-1}:n} | \text{$\tau_{j-1}$, $k$})} = \nonumber \\ 
    &= \frac{P(\tau_{j-1}+1, \tau_j) \; Q_j^{(k)}(\tau_j + 1) \; \pi_k(\tau_{j-1}| \tau_j)}{Q_{j-1}^{(k)}(\tau_{j-1})},
\end{align}
where
 \begin{align}   
    \mathbf{Pr}&(\tau_1 | y_{1:n}, k ) = \label{eq:posteriour_tau_1} \frac{P(1, \tau_1) \; Q^{(k)}(\tau_1 + 1) \; \pi_k(\tau_{1})}{Q^{(k)}(1)}.
\end{align}
Here, $\pi_k$ is the probability of $\tau_j$ based on the distance to $\tau_{j-1}$. This posterior distribution indicates the probability of change point $\tau_j$ occurring in each possible time step $t_i \in [1, n-1]$.
The formulas in \eqref{eq:posteriour_tau_j} and \eqref{eq:posteriour_tau_1} can be applied for each possible number of change points, where $k$ can range from $1$ to $n-1$.
Therefore, this posterior distribution is calculated for
every available number of change points $k$.
\\
\\
The final step in the Bayesian approach is to combine the conditional probabilities for each individual change point (seen in equation \eqref{eq:posteriour_tau_j}) to get the joint distribution for all available change points. The joint probability is calculated as
\begin{equation}
    \mathbf{Pr}(\tau_1, \tau_2, \dotsc, \tau_{n-1} | y_{1:n} ) = \Big( \prod_{j=2}^{n-1} \mathbf{Pr}(\tau_j | \tau_{j-1}, y_{1:n}, k ) \Big) \; \mathbf{Pr}(\tau_1 | y_{1:n}, k),
    \label{eq:posterior_joint}
\end{equation}
where the first change point $\tau_1$ has a different probability formulation due to not having any previous change point. We can also note that the product is changed to a sum if logarithmic probabilities are used, as in \eqref{MAP}. This joint probability can be used to identify the most likely change points. Examples of calculated posterior distributions are found in Appendix A \cite{masterThesis}, where we see the varying probability of being a change point for each sample in the dataset. A sampling method can be used to draw samples from the joint posterior distribution, where we are interested in the points that are most likely to be change points. This means that we can identify the peaks in the posterior distribution, above a set confidence level. This is explained further in section \ref{sec:testing_procedure}.

\subsection{Methods of error estimation}
\label{sec:test_metrics}
In this section, the used metrics for evaluating the performance of the CPD algorithms are presented. We first differentiate between the true change points and the estimated ones. The true change points are denoted by $\mathcal{T}^* = \{\tau^*_0, \dots , \tau^*_{K+1} \}$ while $\Hat{\mathcal{T}}= \{\Hat{\tau}_0, \dots , \Hat{\tau}_{K+1} \}$ indicate estimations. Similarly, the number of true change points is indicated $K^*$ while $\hat{K}$ represents the number of predicted points. 
\\
\\
The most straight forward measure is to compare the number of predictions with the true number of change points. This is know as the \textit{Annotation error}, and is defined as
\begin{equation}
    AE := | \hat{K} - K^* |,
    \label{eq:annotation_error}
\end{equation}
where $\hat{K}$ is the estimated and $K^*$ the true change points. This does not indicate how precise the estimations are, but can indicate if the model is over- or under-fitted.
\\
\\
Another similarity metric of interest is the \textit{Rand Index} (RI) \cite{selective_review_CPD}. Compared to the previous distance metrics, the rand index gives the similarity between two segmentations as a percentage of agreement. This metric is commonly used to compare clustering algorithms.  To calculate the index, we need to define two additional sets which indicate whether two samples are grouped together by a given segmentation or if they are not grouped together. These sets are defined by Truonga et al \cite{selective_review_CPD} as
\begin{align*}
    GR(\mathcal{T}) &:= \{ (s,t), \; 1\leq s < t \leq T : \text{ $s$ and $t$ belong to the same segment in $\mathcal{T}$}\}, \\
    NGR(\mathcal{T}) &:= \{ (s,t), \; 1\leq s < t \leq T : \text{ $s$ and $t$ belong to different segments in $\mathcal{T}$}\},
\end{align*}
where $\mathcal{T}$ is some segmentation for a time interval $[1,T]$. Using these definitions, the rand index is calculated as
\begin{equation*}
    RI(\hat{\mathcal{T}}, \mathcal{T}^*) := \frac{|GR(\hat{\mathcal{T}}) \cap GR(\mathcal{T}^*) | + |NGR(\hat{\mathcal{T}}) \cap NGR(\mathcal{T}^*)|}{T(T-1)}
    \label{eq:randindex}
\end{equation*}
which gives the number of agreements divided by possible combinations. 
\\
\\
To better understand how well the predictions match the actual change points, one can use the measure called the \textit{meantime error} which calculates the meantime between each prediction to the closest actual change point. The meantime should also be considered jointly with the dataset because the same magnitude of meantime error can indicate different things in different datasets. For real-life time series data, the meantime error should be recorded in units of time, such as seconds, in order to make the results intuitive for the user to interpret. The meantime is calculated as
\begin{equation*}
    MT(\hat{\mathcal{T}}, \mathcal{T}^*) = \frac{ \sum_{j = 1}^{\hat{K}} \min_{\tau^* \in \mathcal{T}^*} |\hat{\tau_j} - \tau^* | }{\hat{K}}.
    \label{eq:meantime}
\end{equation*}
A drawback with this measure is that it focuses on the predicted points. If there are fewer predictions than actual change points, the meantime might be lower if the predictions are in proximity of some of the actual change points but not all. Note that the meantime is calculated from the prediction and does not necessarily map the prediction to corresponding true change point, only the closest one. 
\\
\\
Two of the most common metrics of accuracy in predictions are \textit{precision} and \textit{recall}. These metrics give a percentage of how well the predictions reflect the true values. The precision metric is the fraction of correctly identified predictions over the total number of predictions, while the recall metric compares the number of identified true change points over the total number of true change points. These metrics can be expressed as
\begin{equation}
    {\displaystyle {precision}={\frac {|\text{TP}(\hat{\mathcal{T}}, \mathcal{T}^*)|}{|\hat{\mathcal{T}|}}}},
    \quad \quad \quad
    {\displaystyle {recall}={\frac {|\text{TP}(\hat{\mathcal{T}}, \mathcal{T}^*)|}{|\mathcal{T}^*|}}},
    \label{eq:perc_recall}
\end{equation}
where TP represents the number of true positives between the estimations $\hat{\mathcal{T}}$ and true change points $\mathcal{T}^*$. Mathematically, TP is defined as
\begin{equation*}
    \text{TP}(\hat{\mathcal{T}}, \mathcal{T}^*) = \{ \tau^* \in \mathcal{T}^* | \hat{\tau} \in \hat{\mathcal{T}}: |\tau^* - \hat{\tau}|< \epsilon \},
\end{equation*}
where $\epsilon$ is some chosen threshold. The threshold gives the radius of acceptance, meaning the acceptable number of time steps which can differ between prediction and true value. The two metrics \eqref{eq:perc_recall} can be incorporated into a combined metric, known as the \textit{F-score}.
The metric \textit{F1-score} uses the harmonic mean of the precision and recall and is  applied in this work.
As reviewed in this section, the metrics measure the similarity between the predicted change points and the actual change points from various perspectives. Hence this work adopt all of them to give a comprehensive evaluation of the performance of CPD algorithms.

\section{Methods}
\label{sec:methods}
In this section we explore the setting in which the tests are
preformed, along with the testing procedure. First, a description of
the simulated datasets is provided, along with mathematical formulas
and assumptions. Then, we explore the real world dataset with four
process variables. Finally, the testing procedure is described along
with adjustments made for a fair comparison or to reduce computational
complexity. All datasets are described mathematically and illustrated
in figures with the true change points indicated as the boarder
between two segments. All tests shown in this work can be reproduced
using the \texttt{GIT} repositories presented in Appendix B in \cite{masterThesis}.

\subsection{Simulation of data}
\label{sec:simulation_of_data}
To investigate the performance of the approaches with certain features present in the data, simulated datasets might be beneficial to use. The complexity of the datasets can vary and this work studies six simulated datasets. The first four datasets investigate the performance in piecewise constant, piecewise linear, changing variance and autoregressive data respectively. The fifth and sixth datasets indicate realistic processes, with periodic phenomena and non-linear behaviours. Each dataset is explained individually in the following sections. 

\subsubsection*{Piecewise constant}
\begin{figure}[H]
    \hspace{-0.8cm}
    \includegraphics[width = 1.05\linewidth]{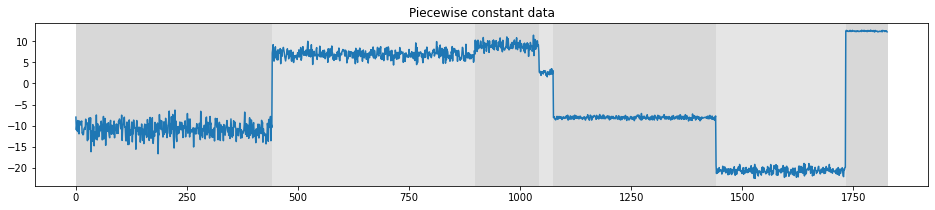}
    \caption{Dataset with seven independent segments, with randomised mean and variance. The segments are indicated with alternating grey backgrounds, and six change points are present on the boarder between segments. }
    \label{fig:piecewise_data}
\end{figure}

To generate the simulated data, we have
created segments with randomised traits (namely mean and variance) and concatenate to get a segmented dataset. If we randomise a mean and variance, we can create a piecewise constant dataset; an example of such data is shown in Figure~\ref{fig:piecewise_data}. In this dataset, we have changes in the mean and variance occurring simultaneously, meaning the mean and the variance of each segment are different from the mean and the variance of other segments. Each value $y_t$ in segment $S_j$ follows the Gaussian distribution
\begin{equation*}
    y_t \sim \mathcal{N}(\mu_j, \sigma_j), \quad j \in \{1, ... , K, K+1\},
\end{equation*}
where $\mu_j \sim \mathcal{U}(-10, 10)$ and $\sigma_j \sim \mathcal{U}(-1,1)$ are randomised constants for each segment. This dataset should be possible to use for computation of CPD in both optimisation and Bayesian approaches, as well as for all cost functions in the optimisation approach.
This dataset may be one of the most manageable datasets. 

\subsubsection*{Piecewise linear}
\begin{figure}[H]
    \hspace{-0.8cm}
    \includegraphics[width = 1.05\linewidth]{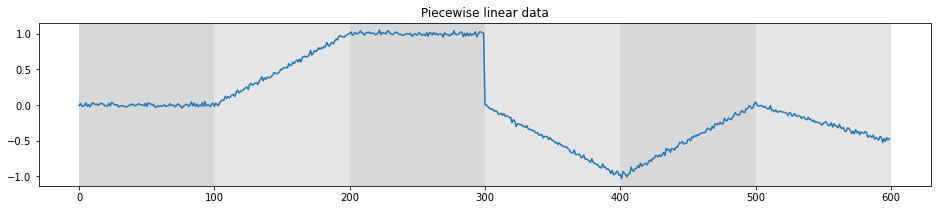}
    \caption{A dataset with piecewise linear segments, where mean is not necessarily constant, and could instead have a linear increase or decrease. The entire dataset has added constant noise. The segments are indicated with alternating grey backgrounds, and five change points are present on the boarder between segments. }
    \label{fig:lin_data}
\end{figure}

In processes, linear changes are common in datasets, when levels transition from one value to another. An example is seen in Figure~\ref{fig:lin_data}, where we see the transition from one linear slope to another at the change points. The dataset is generated by creating linear segments and adding a noise level to the entire dataset, where only the noise level is drawn from a distribution. The difficulty in these datasets lies in the constantly changing mean values in the slopes. Some of the cost functions in the optimisation approach identify changes in the mean values, which would indicate multiple change points along the slopes. Similarly, the Bayesian approach is also attentive to changes in the mean, and would presumably give indications along the slopes. 

\subsubsection*{Changing variance}
\begin{figure}[H]
    \hspace{-0.6cm}
    \includegraphics[width = 1.04\linewidth]{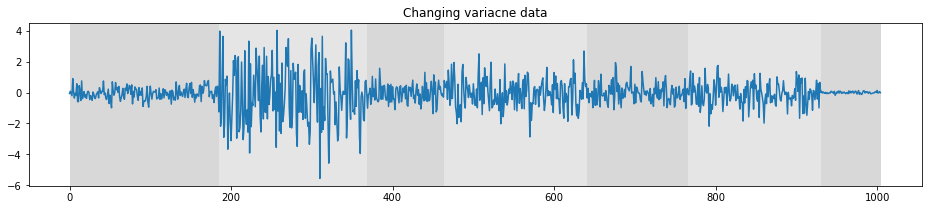}
    \caption{Dataset with constant mean but changing variance. The segments are indicated with alternating grey backgrounds, and six change points are present on the boarder between segments. }
    \label{fig:changing_var_dataset}
\end{figure}

It is interesting to discover the following question:
if the mean is held constant, can detection be performed based on abrupt changes in the variance? The dataset presented in Figure~\ref{fig:changing_var_dataset} is simulated in a similar way to the piecewise constant dataset shown in Figure~\ref{fig:piecewise_data}, except the mean $\mu_j = 0,\; j \in \{1, ..., K+1 \}$. A difficulty in this setting is: how to identify a change point ? Definition~\ref{def:CP} requires presence of a significant change in feature for detection of a change point.
In the setting  of data presented in Figure~\ref{fig:changing_var_dataset}, there are change points indicated without a significant change in variance, which could make CPD more complex. 

\subsubsection*{Autoregressive data}
\begin{figure}[H]
    \hspace{-0.8cm}
    \includegraphics[width = 1.05\linewidth]{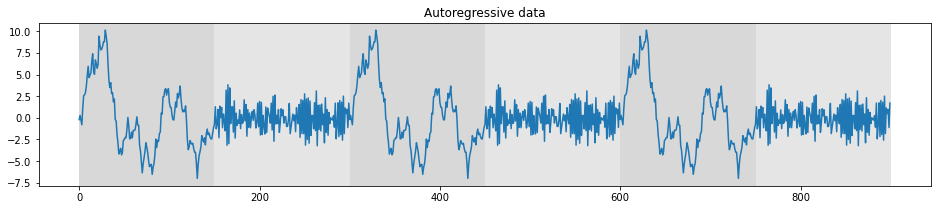}
    \caption{Autoregressive dataset, with six segments and six change points. The segments are indicated with alternating grey backgrounds, and five change points are present on the boarder between segments. }
    \label{fig:AR_dataset}
\end{figure}

The cost functions presented in equations \eqref{eq:c_l1}-\eqref{eq:c_AR} rely on the assumption of an underlying linear model distribution in  data. Especially equation \eqref{eq:c_AR} relies on previous samples and lags, enabling it to fit \textit{autoregressive} (AR) models to segments. This sheds light on the possibility of identifying autoregressive segments and hence change points in AR data. The dataset shown in Figure~\ref{fig:AR_dataset} is constructed using the AR-model with two coefficient (one lag and one constant term). A sample $y_t$ in segment $S_j$ is generated as
\begin{equation*}
    y_{t}= c_j + \varphi^j _{1} y_{{t-1}}+\varepsilon_{t},
\end{equation*}
where $c_j$ is a constant and $\varepsilon_t \sim \mathcal{N}(0,1)$ is the noise. The coefficient $\varphi^j _{1} \in \mathbb{R}$ is a model specific parameter.
Figure~\ref{fig:AR_dataset} presents data generated via the equation above, which are repeated three times for different time segments.
As this dataset also includes a clear change in mean and variance, the performance of the different cost functions are presumed to vary. The Bayesian approach is also expected to have difficulty in identifying specific change points.  

\subsubsection*{Exponential decay data}
\begin{figure}[H]
    \hspace{-0.8cm}
    \includegraphics[width = 1.05\linewidth]{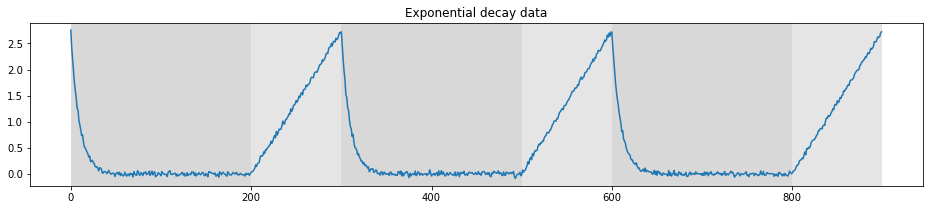}
    \caption{Dataset with exponential decay, followed by a linear increase. Segments are indicated with alternating grey backgrounds. In total there are three repeated processes, with six segments and six change points, each at a boarder between indicated segments.}
    \label{fig:exp_lin_dataset}
\end{figure}

 Now we study 
more realistic features, when some signals
 can mimick  data from chemical processes. It is common to have a segment of exponential decay and a linear segment as a representation of some part of the process. Such behaviour can be of practical relevance. For example, the concentration of a chemical in a reactor can increase linearly when the feed flow of this chemical enters the reactor. Then when the reaction starts, the concentration of this chemical decays exponentially. The change points between these segments indicate the start and the end of the feed flow injection phase and the reaction phase. An illustration of such a process is shown in Figure~\ref{fig:exp_lin_dataset} where three phases are seen constituting of an exponential decay followed by a linear increase. Similar to the piecewise linear dataset, the signal is created and then noise is added. 

\subsubsection*{Oscillating dataset}
\begin{figure}[H]
    \hspace{-0.8cm}
    \includegraphics[width = 1.05\linewidth]{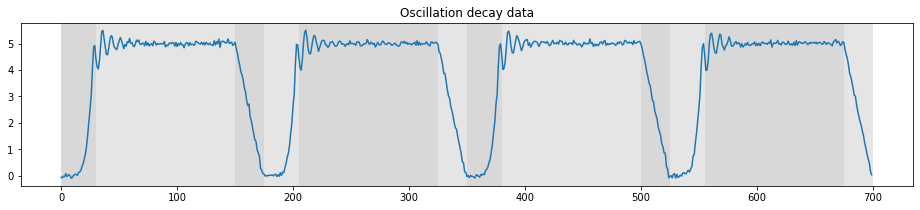}
    \caption{Dataset with a sigmoid function with an oscillation occurring around the maximum value level. The damped oscillation represents the stabilisation of a process. After a stable period, we see a linear decay before the process is repeated. The segments are indicated with alternating grey backgrounds, and eleven change points are present on the boarder between segments. }
    \label{fig:oscil_data}
\end{figure}

Another common phenomenon in processes is a stabilising process when a certain level of stabilisation is reached. This can be represented as a damped oscillation
\begin{equation*}
    y_t = e^{-d t} \cdot \cos(t),
\end{equation*}
where $d\geq 0$ is a damping constant.
In real world applications it is interesting to detect  the point where the stable level is reached, but does not indicates a significant change in features and is, therefore, not a true change point according to the Definition \ref{def:CP}. Figure~\ref{fig:oscil_data}  illustrates a scaled sigmoid function with added oscillations when the target level is reached. Such oscillatory and stabilising behaviour can be often  seen in controlled variables in chemical processes. The two features making this dataset more complex are the sigmoid function and the oscillations occurring before stabilisation. 

\subsection{PRONTO data exploration}
\label{sec:PRONTO}
\begin{figure}[h]
    \centering
    \includegraphics[width = 1 \linewidth]{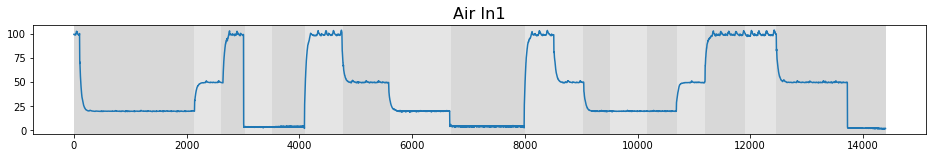}
    \includegraphics[width = 1 \linewidth]{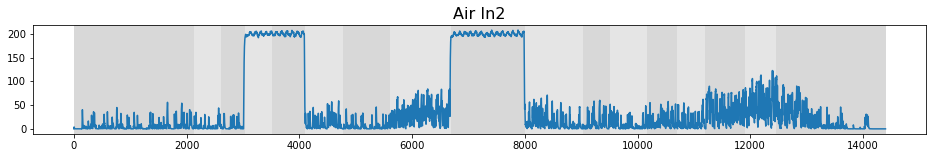}
    \includegraphics[width = 1 \linewidth]{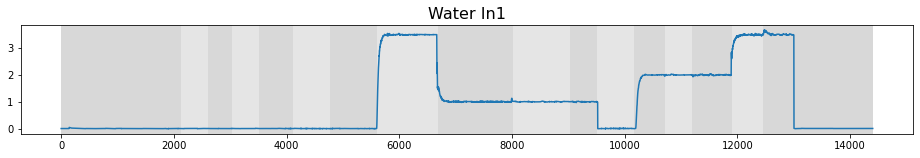}
    \includegraphics[width = 1 \linewidth]{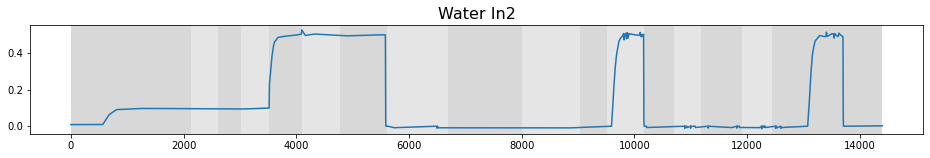}
    \caption{The process variables \texttt{Air In 1}, \texttt{Air In 2}, \texttt{Water In 1} and  \texttt{Water In 2} described in the technical report by PRONTO \cite{PRONTO_description}. The segments between change points are indicated with alternating grey regions, and sixteen change points are present on the boarder between segments.}
    \label{fig:PRONTO_data}
\end{figure}

Multiphase flow processes are frequently used in industries, when two or more substances, such as water, air and oil, are mixed or interact with each other. An example of such a process is described in a technical report conducted within the \textit{Process Network Optimization} project, abbreviated PRONTO \cite{PRONTO_description}. In the described process, air and water are pressurised respectively, where the pressurised mix travels upwards to a separator located at an altitude. Five experimental scenarios are conducted in the facility to monitor the reaction in various process areas, where one experiment is changing the pressure to one of the components. If liquid builds up at the bottom of the facility, this will block the gas flow. The blocking of gas will be  presented until the pressure of the gas is high enough to push the built up liquid to the top. This phenomenon is known as slugging and is an intermittent fault which results in abnormal behaviour in mainly the air and water pressures. For more details about slugging see Figure 1 in 
\cite{PRONTO_description}  which
gives an overview of the facility and the process flow.
\\
\\
The facility usually contains multiple sensors and monitoring systems,
while this work focuses on a few. The data used in this thesis can be
retrieved at \cite{PRONTOdataset}. As the slugging mainly affects the
flows of air and water, their values are used for change point
detection. For each of the components, there are two sensors which can
 measure the respective flow, where one is mainly used and the
second one is used only in some cases. This means that we will focus
on four process variables, two for each component. The process
variables are denoted  by \texttt{Air In 1}, \texttt{Air In
2}, \texttt{Water In 1} and \texttt{Water In
2}, respectively, see Figure ~\ref{fig:PRONTO_data}.  Figure~\ref{fig:PRONTO_data} shows these four signals, where the
segments between change points are indicated with alternating grey
colour. We observe a range of the features presented in the simulated
datasets of all signals such that piecewise constant segments,
exponential decay and change in variance. The signals are sampled with
the same sampling rate and can be examined simultaneously, but will be
treated individually in this work. This means that predictions are
made for each signal individually and these predictions are then
aggregated, namely the final detected change points are the union of
all change points in all signals.  This is done to compare the
predictions for the actual change points, which are not necessarily
linked for only one process variable.

\subsection{Testing procedure}
\label{sec:testing_procedure}
To perform accurate testing on the previously described datasets, a detailed procedure is needed. The two studied approaches - optimisation and Bayesian - have different procedures to predict change points, and the procedure of identifying the best predictions varies in the two approaches. Below, the two procedures as well
as the set-up of parameters for tested problems are explained in detail.
\\
\\
To make predictions with the \texttt{RUPTURES} package, first one needs to specify the algorithm with search direction and cost function. In addition to this, the penalty level should also to be specified. The value of the perfect penalty $\beta$ in equation \eqref{eq:opt_problem} is not known a priori, therefore, multiple predictions are done for various penalty values $\beta \in \mathbb{N} \cup \{ 0 \}$. Note, if $\beta = 0$ the algorithm receives no penalty when adding a change point, and the algorithm can add as many change points as necessary to minimise the cost functions. Predictions are made using all cost functions, where both \texttt{PELT} and \texttt{WIN} are used, and the resulting predictions are saved. For each prediction, the metrics provided in section \ref{sec:test_metrics} are calculated and saved. When enough tests have been performed in terms of penalty values, the results are saved to an external file. To select the best prediction, the metric values need to be taken into account. Our goal is
to minimise the annotation error and meantime, and
at the same time
maximise the F1-score and the rand index. Discussion on how to evaluate the metrics and choose the best one for different data is provided in section \ref{sec:discussion}. Obtained results are presented in  section \ref{sec:results} along with respective metric values. 
\\
\\
To apply the Bayesian approach for predictions, the  computational procedure is different.
Using the concepts derived in section \ref{sec:bayes} we can calculate posterior distribution with probabilities for each time step being a change point. Due to the algorithm being computationally heavy, the resolution of the data is reduced in the real dataset by PRONTO using \textit{Piecewise Aggregate Approximation} (PAA) \cite{PAA}. The function aggregates the values of a window to an average value. The used window size is $20$ samples. Generally, to draw conclusions from a posterior distribution, sampling is used to create a collection of points which in this case represents the change points. The calculated distribution does not follow a simple distribution, which makes sampling complicated. In essence, we want to create a sample of the most probable change points, without unnecessary duplicates. To draw this type of sample\footnote{This is not necessarily a proper sampling methodology, and other approaches can be used instead. An alternative sampling method is provided in Fearnhead \cite{Exact_efficient_Bayesian_inference} (page 8). } of change points from the posterior distribution, the function \texttt{find\_peaks} in the \texttt{Python} package \texttt{SciPy}\cite{find_peaks} is used. The function identifies the peaks in a dataset using two parameters: \texttt{threshold} which the peak value should exceed, and \texttt{distance} which indicates the minimum distance between peaks. The threshold is set to $0.2$, where we require a certainty level of at least $20\%$. The distance is set to $10$ time steps to prevent duplicate values. The posterior distribution is calculated once for one dataset, where numerous samples can be drawn using different settings in the \texttt{find\_peaks} function. This approach returns the most probable change points which are then used to calculate the metrics presented in section~\ref{sec:test_metrics}. 
\\
\\
All signals are handled individually meaning we are only investigating the uni-variate case, without correlation between the covariates. In the simulated datasets, this is trivial since we only have one signal per case. In the PRONTO dataset we have four process variables, which are explained in the previous section. The same prediction algorithm is applied to all signals, and are not altered between the different process variables. This means that the range of the signals can affect the predictions. To counteract unfair predictions, the process variables are normalised. Normalisation is not necessary for the signals in the simulated datasets, while the process variables in the PRONTO dataset are normalised to account for the difference in range in the signals.

\section{Results}
\label{sec:results}
Given the different search directions and cost functions, presented in sections \ref{sec:opt_search_dir} and \ref{sec:opt_cost_functions} respectively, we can presume that different setups will identify different features and hence differ in the prediction of change points. We can also assume that the Bayesian approach, presented in section \ref{sec:bayes}, will not necessarily give the same predictions as the optimisation approach. We note that all algorithms predict the intermediate change points $\mathcal{T} = \{\tau_1,...,\tau_K\}$ along with one artificial change point $\tau_{K+1} := n$. This artificial change point is based on definition and is used when the predictions are compared. A first step to understanding the performance of the different approaches is to simulate datasets with certain features and compare the obtained metrics. In addition, a visualisation is shown for each case, with the predicted change point in comparison to the actual change points. In this section, we present the results of the two approaches on the six simulated datasets with varying complexity, described in section \ref{sec:simulation_of_data}.
Later, the results for the real-world data are presented. 

\subsection{Simulated datasets}
The results for the simulated datasets are presented individually along with an illustration of the dataset and predictions.
Results for both approaches are presented for each dataset separately in Tables 1, 2 and 3. Table~\ref{tab:res_WIN_simulated_datasets}~and~\ref{tab:res_PELT_simulated_datasets} present the predictions made using the approximate and exact search directions, respectively. The results from the Bayesian approach can be found in Table~\ref{tab:res_BAYES_simulated_datasets}.
\\
\\
\begin{table}[H]
    \centering
    \begin{tabular}{l| c c | c c c c c}
       Dataset & K & AE & MT & Precision & Recall & F1 & RI \\
        & && [time steps] & [\%] & [\%] & [\%] & [\%] \\
        \hline
        Piecewise constant & 7 & 0 & 9.8 & 83.3 & 83.3 & 83.3 & 98.5 \\
        
        Piecewise linear & 2 & 4 & 2.0 & 100 & 20.0 & 33.3 & 77.5 \\
        
        Changing variance & 4 & 3 & 2.0 & 100 & 50.0 & 66.7 & 76.5 \\
        
        Autoregressive data & 19 & 13 & 46.7 & 27.8 & 100 & 43.5 & 94.5 \\
        
        Exponential decay & 9 & 3 & 31.4 & 62.5 & 100 & 76.9 & 92.8 \\
        
        Oscillating data & 12 & 0 & 3.5 & 72.7 & 72.7 & 72.7 & 98.3 \\
        \hline
        
    \end{tabular}
    \caption{Prediction results for the simulated datasets, using the Bayesian method. The indicated change points are the peaks in the predicted posterior probability distribution. }
    \label{tab:res_BAYES_simulated_datasets}
\end{table}

\newpage
\thispagestyle{empty}

\begin{table}[H]
    \vspace{-0.2cm}
    \hspace{-0.8cm}
    \begin{tabular}{l|l c | c c | c c c c c}
        Dataset & Cost & Pen & K & AE & MT & Precision & Recall & F1 & RI \\
        & function &&&& [time steps] & [\%] & [\%] & [\%] & [\%] \\
        \hline
        \multirow{7}{*}{Piecewise constant} & $c_{L2}$ & 5 & 6 & 1  & 1.2 & 100 & 83.3 & 90.9 & 99.5 \\ 
        & $c_{L1}$ & 1 & 7 & 0 & 20.8 & 83.3 & 83.3 & 83.3 & 97.2 \\
        & $c_{Normal}$ & 5 & 6 & 1 & 1.4 & 100 & 83.3 & 90.9 & 99.5 \\ 
        & $c_{LinReg}$ & 424 & 6 & 1 & 57.4 & 80.0 & 66.7 & 72.7 & 92.4 \\
        & $c_{AR}$ & 5 & 5 & 2 & 1.5 & 100 & 66.7 & 80.0 & 94.7 \\
        & $c_{ridge}$ & 24 & 6 & 1 & 11.4 & 100 & 83.3 & 90.9 & 98.1 \\ 
        & $c_{lasso}$ & 24 & 6 & 1 & 11.4 & 100 & 83.3 & 90.9 & 98.1 \\ 
        \hline
        \multirow{7}{*}{Piecewise linear} & $c_{L2}$ & 0 & 5 & 1 & 36.3 & 25.0 & 20.0 & 22.2 & 88.4 \\
        & $c_{L1}$ & 0 & 5 & 1 & 33.8 & 25.0 & 20.0 & 22.2 & 88.7 \\
        & $c_{Normal}$ & 0 & 5 & 1 & 7.5 & 100 & 80.0 & 88.9 & 92.8 \\
        & $c_{LinReg}$ & 0 & 5 & 1 & 35.9 & 25.0 & 20.0 & 22.2 & 88.2 \\
        & $c_{AR}$ & 0 & 6 & 0 & 5.0 & 100 & 100 & 100 & 97.6  \\
        & $c_{ridge}$ & 0 & 6 & 0 & 0.0 & 100 & 100 & 100 & 100 \\
        & $c_{lasso}$ & 0 & 6 & 0 & 1.0 & 100 & 100 & 100 & 100 \\
        \hline
        \multirow{7}{*}{Changing variance} & $c_{L2}$ & 0 & 7 & 0 & 42.2 & 50.0 & 50.0 & 50.0 & 91.3 \\
        & $c_{L1}$ & 0 & 7 & 0 & 45.2 & 66.7 & 66.7 & 66.7 & 90.4 \\
        & $c_{Normal}$ & 0 & 7 & 0 & 22.0 & 83.3 & 83.3 & 83.3 & 92.0 \\
        & $c_{LinReg}$ & 0 & 7 & 0 & 42.2 & 50.0 & 50.0 & 50.0 & 91.3 \\
        & $c_{AR}$ & 0 & 6 & 1 & 44.8 & 60.0 & 50.0 & 55.0 & 87.4 \\
        & $c_{ridge}$ & 0 & 3 & 4 & 40.5 & 50.0 & 16.7 & 25.0 & 74.9 \\
        & $c_{lasso}$ & 0 & 3 & 4 & 40.5 & 50.0 & 16.7 & 25.0 & 74.9 \\
        \hline
        \multirow{7}{*}{Autoregressive data} & $c_{L2}$ & 0 & 7 & 1 & 37.5 & 83.3 & 100 & 90.9 & 91.3 \\
        & $c_{L1}$ & 0 & 7 & 1 & 40.0 & 83.3 & 100 & 90.9 & 90.5 \\
        & $c_{Normal}$ & 0 & 7 & 1 & 40.0 & 83.3 & 100 & 90.9 & 90.5 \\
        & $c_{LinReg}$ & 0 & 7 & 1 & 37.5 & 83.3 & 100 & 90.9 & 91.3 \\
        & $c_{AR}$ & 0 & 6 & 0 & 2.0 & 100 & 100 & 100 & 99.3 \\
        & $c_{ridge}$ & 0 & 4 & 2 & 71.7 & 0.0 & 0.0 & 0.0 & 79.4 \\
        & $c_{lasso}$ & 0 & 4 & 2 & 71.7 & 0.0 & 0.0 & 0.0 & 79.4 \\
        \hline
        \multirow{7}{*}{Exponential decay} & $c_{L2}$ & 0 & 4 & 2 & 21.7 & 66.7 & 40.0 & 50.0 & 89.7 \\
        & $c_{L1}$ & 0 & 5 & 1 & 27.5 & 50.0 & 40.0 & 44.4 & 89.9 \\
        & $c_{Normal}$ & 0 & 6 & 0 & 16.0 & 100 & 100 & 100 & 95.1 \\
        & $c_{LinReg}$ & 0 & 6 & 0 & 30.0 & 60.0 & 60.0 & 60.0 & 92.5 \\ 
        & $c_{AR}$ & 0 & 7 & 1 & 18.3 & 83.3 & 100 & 90.9 & 98.5 \\
        & $c_{ridge}$ & 0 & 6 & 0 & 6.0 & 100 & 100 & 100 & 97.9 \\
        & $c_{lasso}$ & 0 & 6 & 0 & 36.0 & 80.0 & 80.0 & 80.0 & 89.8 \\
        \hline
        \multirow{7}{*}{Oscillating data} & $c_{L2}$ & 0 & 5 & 7 & 7.5 & 100 & 36.4 & 53.3 & 87.3 \\
        & $c_{L1}$ & 0 & 5 & 7 & 8.8 & 100 & 36.4 & 53.3 & 88.3 \\
        & $c_{Normal}$ & 0 & 8 & 4 & 12.9 & 100 & 63.6 & 77.8 & 93.0 \\
        & $c_{LinReg}$ & 0 & 4 & 8 & 0.0 & 100 & 27.3 & 42.9 & 86.3 \\ 
        & $c_{AR}$ & 0 & 5 & 7 & 10.0 & 100 & 36.4 & 53.3 & 88.9 \\
        & $c_{ridge}$ & 0 & 5 & 7 & 8.8 & 100 & 36.4 & 53.3 & 88.9 \\
        & $c_{lasso}$ & 0 & 5 & 7 & 8.8 & 100 & 36.4 & 53.3 & 88.9\\ 
        \hline
        
    \end{tabular}
    \caption{Prediction results for the optimisation approach, using the search method \texttt{WIN}. For each cost function, a possible penalty level is indicated to obtain the prediction. The used window width $w = 100$ is for all approaches. }
    \label{tab:res_WIN_simulated_datasets}
\end{table}

\begin{table}[H]
    \vspace{-0.2cm}
    \hspace{-0.8cm}
    \begin{tabular}{l|l c | c c | c c c c c}
       Dataset & Cost & Pen & K & AE & MT & Precision & Recall & F1 & RI \\
        & function &&&& [time steps] & [\%] & [\%] & [\%] & [\%] \\
        \hline
        \multirow{7}{*}{Piecewise constant} & $c_{L2}$ & 450 & 7 & 0 & 1.7 & 66.7 & 66.7 & 66.7 & 95.8 \\
        & $c_{L1}$ & 15 & 7 & 0 & 1.3 & 83.3 & 83.3 & 83.3 & 99.8 \\
        & $c_{Normal}$ & 10 & 6 & 1 & 3.4 & 100 & 83.3 & 90.9 & 95.0 \\
        & $c_{LinReg}$ & 1000 & 7 & 0 & 99.2 & 66.7 & 66.7 & 66.7 & 92.6 \\
        & $c_{AR}$ & 80 & 7 & 0 & 1.8 & 66.7 & 66.7 & 66.7 & 94.8 \\
        & $c_{Ridge}$ & 5 & 7 & 0 & 18.8 & 83.3 & 83.3 & 83.3 & 98.3 \\
        & $c_{Lasso}$ & 70 & 7 & 0 & 1.3 & 83.3 & 83.3 & 83.3 & 99.8 \\
        \hline
        \multirow{7}{*}{Piecewise linear} & $c_{L2}$ & 2 & 6 & 0 & 31.0 & 40.0 & 40.0 & 40.0 & 89.4 \\
        & $c_{L1}$ & 6 & 6 & 0 & 30.0 & 60.0 & 60.0 & 60.0 & 89.8 \\
        & $c_{Normal}$ & 120 & 6 & 0 & 23.0 & 60.0 & 60.0 & 60.0 & 91.8 \\
        & $c_{LinReg}$ & 0 & 5 & 1 & 10.0 & 75.0 & 60.0 & 66.7 & 93.1 \\
        & $c_{AR}$ & 0 & 6 & 0 & 4.0 & 80.0 & 80.0 & 80.0 & 93.1 \\
        & $c_{Ridge}$ & 0 & 6 & 0 & 0.0 & 100 & 100 & 100 & 100 \\
        & $c_{Lasso}$ & 2 & 6 & 0 & 1.0 & 100 & 100 & 100 & 99.5 \\
        \hline
        \multirow{7}{*}{Changing variance} & $c_{L2}$ & 6 & 8 & 1 & 46.4 & 28.6 & 33.3 & 30.8 & 63.5 \\
        & $c_{L1}$ & 3 & 7 & 0 & 47.3 & 16.7 & 16.7 & 16.7 & 58.8 \\
        & $c_{Normal}$ & 9 & 6 & 1 & 3.4 & 100 & 83.3 & 90.9 & 95.0 \\
        & $c_{LinReg}$ & 0 & 8 & 1 & 40.7 & 57.1 & 66.7 & 61.5 & 91.3 \\
        & $c_{AR}$ & 30 & 6 & 1 & 51.2 & 20.0 & 16.7 & 18.2 & 59.3 \\
        & $c_{Ridge}$ & 0 & 9 & 2 & 41.0 & 62.5 & 83.3 & 71.4 & 90.7 \\
        & $c_{Lasso}$ & 2 & 8 & 1 & 45.6 & 28.6 & 33.3 & 30.8 & 69.5 \\
        \hline
        \multirow{7}{*}{Autoregressive data} & $c_{L2}$ & 320 &6 & 0 & 40.0 & 40.0 & 40.0 & 40.0 & 83.4 \\
        & $c_{L1}$ & 60 & 6 & 0 & 40.0 & 40.0 & 40.0 & 40.0 & 83.4 \\
        & $c_{Normal}$ & 100 & 6 & 0 & 40.0 & 40.0 & 40.0 & 40.0 & 83.4 \\
        & $c_{LinReg}$ & 145 & 6 & 0 & 28.0 & 80.0 & 80.0 & 80.0 & 93.4 \\
        & $c_{AR}$ & 30 & 6 & 0 & 2.0 & 100 & 100 & 100 & 99.3 \\
        & $c_{Ridge}$ & 10 & 6 & 0 & 50.0 & 40.0 & 40.0 & 40.0 & 86.0 \\
        & $c_{Lasso}$ & 50 & 6 & 0 & 39.0 & 40.0 & 40.0 & 40.0 & 82.9 \\
        \hline
        \multirow{7}{*}{Exponential decay} & $c_{L2}$  & 20 & 6 & 0 & 26.0 & 100 & 100 & 100 & 93.5 \\
        & $c_{L1}$ & 20 & 6 & 0 & 25.0 & 100 & 100 & 100 & 93.6 \\
        & $c_{Normal}$ & 200 & 7 & 1 & 40.0 & 83.3 & 100 & 90.9 & 91.8 \\
        & $c_{LinReg}$ & 5 & 6 & 0 & 12.0 & 100 & 100 & 100 & 96.7 \\
        & $c_{AR}$ & 0 & 8 & 2 & 1.4 & 71.4 & 100 & 83.3 & 99.5 \\
        & $c_{Ridge}$ & 15 & 6 & 0 & 3.0 & 100 & 100 & 100 & 99.1 \\
        & $c_{Lasso}$ & 35 & 6 & 0 & 4.0 & 100 & 100 & 100 & 98.8 \\
        \hline
        \multirow{7}{*}{Oscillating data} & $c_{L2}$ & 15 & 12 & 0 & 6.8 & 72.7 & 72.7 & 72.7 & 97.0 \\
        & $c_{L1}$ & 10 & 12 & 0 & 5.9 & 72.7 & 72.7 & 72.7 & 97.1 \\
        & $c_{Normal}$ & 75 & 12 & 0 & 8.2 & 81.8 & 81.8 & 81.8 & 95.6 \\
        & $c_{LinReg}$ & 0 & 5 & 7 & 12.5 & 75.0 & 27.3 & 40.0 & 87.5 \\ 
        & $c_{AR}$ & 1 & 12 & 0 & 5.0 & 72.7 & 72.7 & 72.7 & 96.8 \\
        & $c_{Ridge}$ & 0 & 6 & 6 & 24.0 & 60.0 & 27.3 & 37.5 & 89.7 \\
        & $c_{Lasso}$ & 10 & 12 & 0 & 1.8 & 72.7 & 72.7 & 72.7 & 98.9 \\
        \hline
        
    \end{tabular}
    \caption{Prediction results for the optimisation approach, using the search method \texttt{PELT}. For each cost function, a possible penalty level is indicated to obtain the prediction along with respective obtained scores. }
    \label{tab:res_PELT_simulated_datasets}
\end{table}
\thispagestyle{empty}

\subsubsection{Piecewise constant data}
\begin{figure}[h]
    \hspace{-0.7cm}
    \includegraphics[width = 1.05\linewidth]{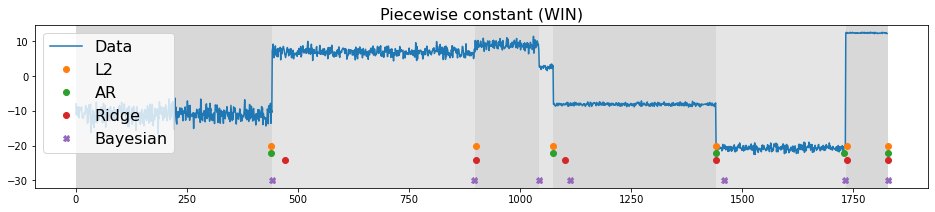}

    \hspace{-0.7cm}
    \includegraphics[width = 1.05\linewidth]{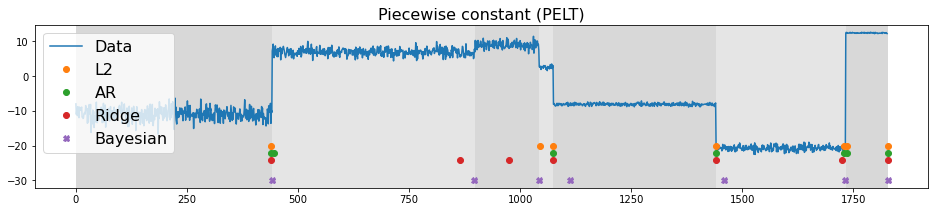}
    \caption{Results for the piecewise constant dataset using various cost functions and the two search directions \texttt{WIN} and \texttt{PELT} respectively. Predicted change points are seen in coloured points, which are time points and do not have a value. The cost functions based on the $L_2$-norm, autoregressive model and Ridge regression can be seen in both illustrations, where the $c_{L2}$ agrees with the actual change points. The predictions made with the Bayesian approach can be seen for comparison to the other predictions.}
    \label{fig:res_piecewise_constant}
\end{figure}
\noindent
As described in section \ref{sec:simulation_of_data}, piecewise constant data is one of the simpler datasets with distinct changes in features. Table~\ref{tab:res_WIN_simulated_datasets} and Table~\ref{tab:res_PELT_simulated_datasets} show that most algorithms predict the change points well, in terms of rand index and F1 score. Using these tables, we observe that the best
predictions are obtained using $L_2$-norm in \texttt{WIN} algorithm,
with a rand index of $99.5\%$ and F1-score of $90.9\%$. Six predictions are made, where one is missed, which in turn affects the recall. The meantime of $1.2$ time steps indicates accurate predictions. We observe also that
the model based cost functions give accurate predictions  as well, where the
AR-model has the smallest meantime with fewest predictions. In Figure~\ref{fig:res_piecewise_constant} (top) we see a comparison between $c_{L2}$, $c_{AR}$, $c_{Ridge}$ and the Bayesian approach, where we observe that the $c_{L2}$ indicates most of the actual change points. Using the \texttt{PELT} search direction, most cost functions identify all seven change points, except $c_{Normal}$ which only identifies six change points. Using $c_{L1}$ and $c_{Lasso}$ the smallest meantime is recorded at $1.3$ timesteps, and the highest RI is $99.8\%$. In Figure~\ref{fig:res_piecewise_constant} (bottom) we see a comparison between $c_{L2}$, $c_{AR}$ and $c_{Ridge}$ which all predict seven change points. Using the $L2$-norm with \texttt{PELT} gives similar predictions as the same cost function with \texttt{WIN}, which agrees with most of the true change points. Cost function $c_{AR}$ has a low meantime, and RI of $94.8\%$, where the algorithm makes two double predictions when using \texttt{PELT}. The ridge regression cost function also predicts seven change points, which deviates from the actual change points in some cases. In both illustrations in Figure~\ref{fig:res_piecewise_constant} the predictions from the Bayesian approach are indicated. The accuracy metrics of the predictions are presented in Table~\ref{tab:res_BAYES_simulated_datasets}, where the Bayesian approach has a meantime of $9.8$, F1-score $83.3\%$ and RI $98.5\%$. 

\subsubsection{Piecewise linear data}
\begin{figure}[h]
    \hspace{-0.7cm}
    \includegraphics[width = 1.05\linewidth]{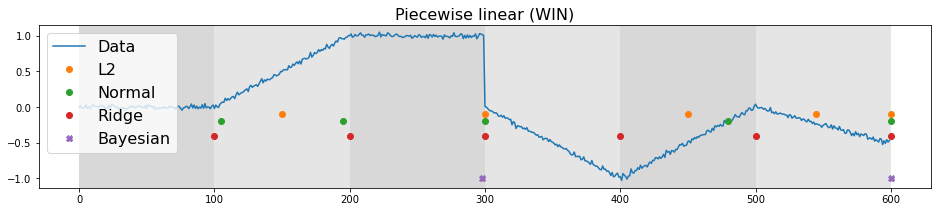}

    \hspace{-0.7cm}
    \includegraphics[width = 1.05\linewidth]{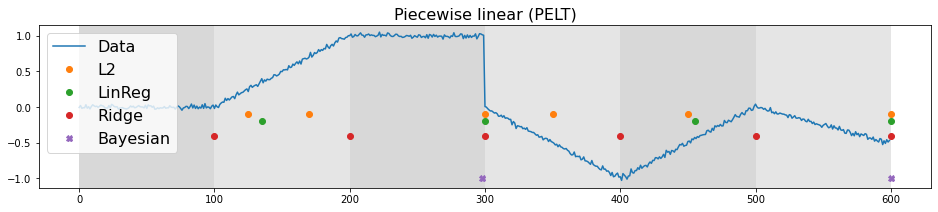}
    \caption{Results for the piecewise constant dataset using various cost functions. Predicted change points are seen in coloured points, which are time points and do not have a value. The cost functions $c_{L2}$, $c_{Normal}$ and $c_{Ridge}$ can be seen in the upper figure, where \texttt{WIN} is used. The bottom image gives a comparison of $c_{L2}$, $c_{Normal}$ and $c_{Ridge}$ when \texttt{PELT} is used. In both illustrations, the ridge regression model agrees well with the actual change points. The predictions made with the Bayesian approach can be seen for comparison to the other predictions.}
    \label{fig:res_piecewise_linear}
\end{figure}
\noindent
Instead of an abrupt change in the mean, piecewise linear data can illustrate how the mean changes continuously, which affects some algorithms' performance. The norm-based cost functions $c_{L2}$ and $c_{L1}$ have the lowest F1-score and RI when \texttt{WIN} is used, along with $c_{LinReg}$ and $c_{Lasso}$, see Table~\ref{tab:res_WIN_simulated_datasets}. The autoregressive cost function $c_{AR}$ predicts all six change points with a meantime of $5$ timesteps, which corresponds to $0.8\%$ of the samples in the data. Similarly, $c_{ridge}$ also predicts all six points, but with a zero meantime error and $100\%$ F1-score. This is also the case for the predictions made with  $c_{ridge}$ and \texttt{PELT}. When the exact search direction is used, some RI and F1-scores are lower than for \texttt{WIN}, except for $c_{L2}$, $c_{L1}$ and $c_{Lasso}$, where predictions are slightly improved in Table~\ref{tab:res_PELT_simulated_datasets}. In Figure~\ref{fig:res_piecewise_linear} we see a comparison of some cost functions with either \texttt{WIN} or \texttt{PELT}. Note that only $c_{Ridge}$ predicts all change points in this case. The Bayesian approach can be seen in the figures, where the approach only predicts two change points. In Table~\ref{tab:res_BAYES_simulated_datasets} we find the F1-score to be $33.3\%$ and the RI $77.5\%$, which is lower than most of the predictions made using either \texttt{WIN} or \texttt{PELT}. 

\subsubsection{Changing variance}
\begin{figure}
    \hspace{-0.7cm}
    \includegraphics[width = 1.05\linewidth]{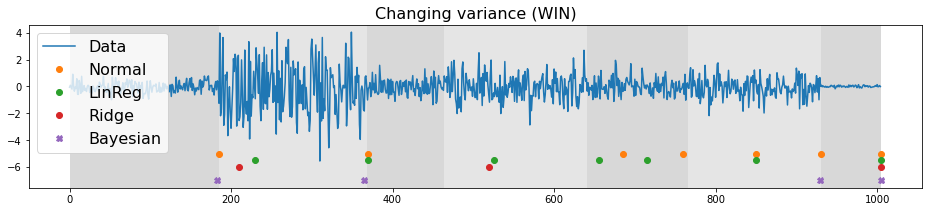}

    \hspace{-0.7cm}
    \includegraphics[width = 1.05\linewidth]{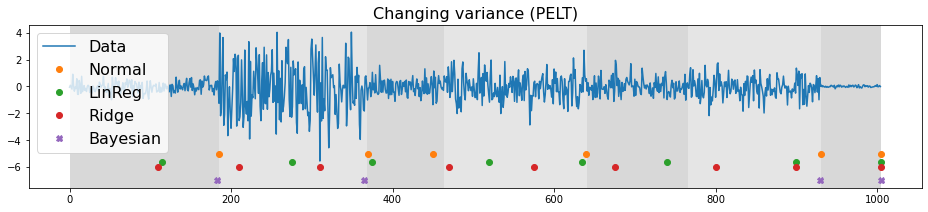}
    \caption{Results for the dataset with constant mean and changing variance, using various cost functions. Predicted change points are seen in coloured points, which are time points and do not have a value. The cost functions $c_{Normal}$, $c_{LinReg}$ and $c_{Ridge}$ can be seen in the upper figure, where \texttt{WIN} is used. The bottom image gives a comparison of $c_{L2}$, $c_{LinReg}$ and $c_{Ridge}$ when \texttt{PELT} is used. All predictions give different predictions compared to the predictions made by the Bayesian approach, where only four change points are predicted. }
    \label{fig:res_changing_var}
\end{figure}
The two previously discussed datasets have features
for following three cases:
when the mean changes throughout the data, when the mean can be held constant and when only the variance changes at change points. In Table~\ref{tab:res_WIN_simulated_datasets} we see how $c_{L2}$, $c_{L1}$, $c_{Normal}$ and $c_{LinReg}$ predict seven change points when using \texttt{WIN}, while $c_{Ridge}$ predicts three change points. Using $c_{Lasso}$ gives no predictions, except the implicit change point at time $T$. We can note that the meantime is generally high and the F1-score is lower than $85\%$ for all cost functions. In Figure~\ref{fig:res_changing_var} (upper) we see the predictions made when using $c_{Normal}$, $c_{LinReg}$ and $c_{Ridge}$, where $c_{Normal}$ agrees with five of the seven change points. When using the exact search direction (\texttt{PELT}) the accuracy in the algorithms' predictions change, as seen in Table~\ref{tab:res_PELT_simulated_datasets}. In this case $c_{LinReg}$ and $c_{Ridge}$ have the highest F1-score of $61.5\%$ and $71.4\%$ respectively. The two algorithms have rand index scores above $90\%$. With \texttt{PELT}, $c_{Normal}$ makes six predictions and F1-score of $90.9\%$ and a rand index of $95\%$. This cost function has the lowest meantime and the highest precision, recall and rand index compared to all other prediction methods. In Figure~\ref{fig:res_changing_var} (bottom), we see the predictions made by $c_{Normal}$, $c_{LinReg}$ and $c_{Ridge}$, where none of the three algorithms agree with the true segmentation. We note the cluster of predictions made by $c_{L2}$, which is different from the other predictions. In both of the images presented in Figure~\ref{fig:res_changing_var} we see the Bayesian predictions, which indicates four change points. Table~\ref{tab:res_BAYES_simulated_datasets} shows a low meantime for the predictions and a high precision. The recall is $50\%$ and in turn results in the F1-score of $66.7\%$.

\subsubsection{Autoregressive data}
\begin{figure}
    \hspace{-0.7cm}
    \includegraphics[width = 1.05\linewidth]{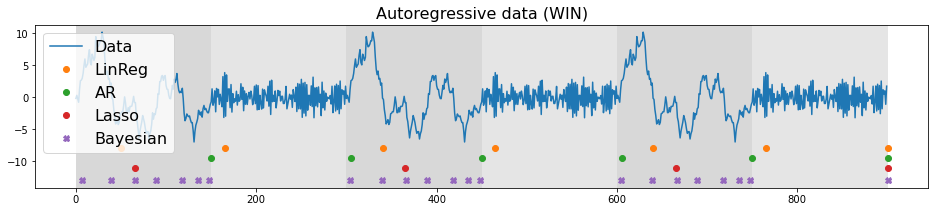}

    \hspace{-0.7cm}
    \includegraphics[width = 1.05\linewidth]{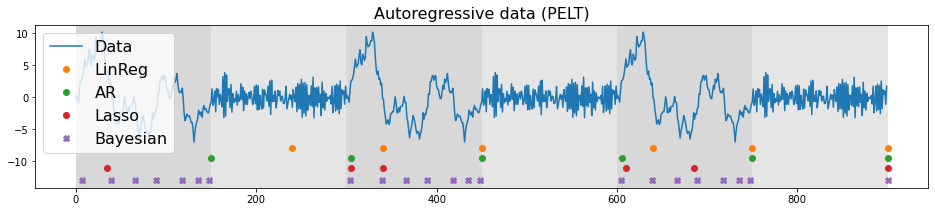}
    \caption{Results for the autoregressive dataset using various cost functions. Predicted change points are seen in coloured points, which are time points and do not have a value. The cost functions $c_{LinReg}$, $c_{AR}$ and $c_{Lasso}$ can be seen in both figures, where $c_{AR}$ agrees well with the true change points. All predictions give different predictions compared to the predictions made by the Bayesian approach, where three clusters of change points are present. }
    \label{fig:res_AR_data}
\end{figure}
In the autoregressive dataset we see a change in variance in the different segments, along with a change on mean depending on which subset is viewed. Using the approximate search direction \texttt{WIN} some cost functions give good predictions of the change points. In Table~\ref{tab:res_WIN_simulated_datasets} we see how $c_{AR}$ predicts six change points with a meantime of two time steps, which is significantly lower than all other metrics. The algorithm has F1-score of $100\%$ and rand index $99.3\%$. The cost functions $c_{L2}$, $c_{L1}$, $c_{Normal}$ and $c_{LinReg}$ predict seven change points, where the meantime is lower for $c_{L2}$ and $c_{LinReg}$. The regularisation functions $c_{Ridge}$ and $c_{Lasso}$ predict the same number of change points, have the lowest F1-scores and rand index as well as the largest meantime. When using \texttt{PELT} the predictions using $c_{AR}$ do not change, while all other cost functions give other predictions. In Table~\ref{tab:res_PELT_simulated_datasets}, all cost functions predict six change points, except $c_{Normal}$ which predicts two. Compared to \texttt{WIN}, the rand index is lower for $c_{L2}$, $c_{L1}$ and $c_{Normal}$, while $c_{LinReg}$, $c_{Ridge}$ and $c_{Lasso}$ show an increase in Table~\ref{tab:res_PELT_simulated_datasets}. Figure~\ref{fig:res_AR_data} we see the predictions of some algorithms using \texttt{WIN} (top) and \texttt{PELT} (bottom) respectively. In both images, we see how $c_{AR}$ agrees with all the true change points, while $c_{LinReg}$ and $c_{Lasso}$ do not. We can also see the predictions made by the Bayesian approach, where we note the clustering of predictions. In Table~\ref{tab:res_BAYES_simulated_datasets} we see that $19$ change points are predicted, with a meantime of $46.7$ time steps. The F1-score is $43.5\%$ and rand index $94.5\%$. 

\subsubsection{Exponential decay data}
\begin{figure}
    \hspace{-0.7cm}
    \includegraphics[width = 1.05\linewidth]{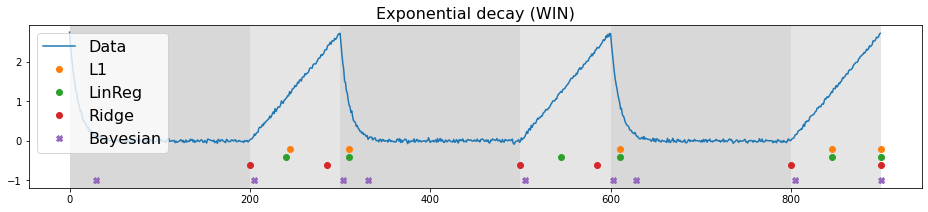}

    \hspace{-0.7cm}
    \includegraphics[width = 1.05\linewidth]{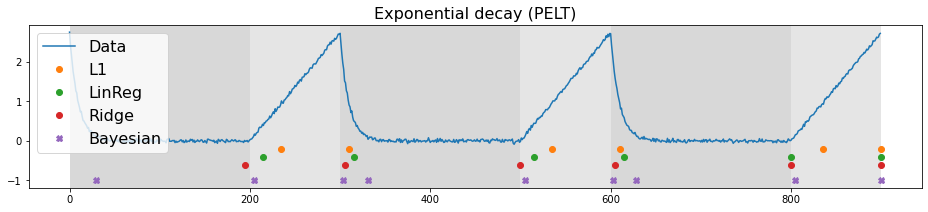}
    \caption{Results for the dataset with exponential trends present, using various cost functions. Predicted change points are seen in coloured points, which are time points and do not have a value. The cost functions $c_{L1}$, $c_{LinReg}$ and $c_{Ridge}$ can be seen in both figures, where $c_{Ridge}$ agrees with the true change points using both \texttt{WIN} and \texttt{PELT}. The predictions resemble the predictions made by the Bayesian approach, where three additional change points are predicted. }
    \label{fig:res_exp_data}
\end{figure}
As an extension to the piecewise linear data, exponential trends can be present in datasets. In Figure~\ref{fig:res_exp_data} we see this piecewise linear and exponential decay in combination with the predictions made by various approaches. The upper figure shows the predictions made using $c_{L1}$, $c_{LinReg}$ and $c_{Ridge}$ using \texttt{WIN}. We see how $c_{Ridge}$ gives a prediction agreeing with the true change points. Table~\ref{tab:res_WIN_simulated_datasets} presents the highest F1-scores obtained by $c_{Normal}$, $c_{Ridge}$ and $c_{Lasso}$, where $c_{Ridge}$ and $c_{Lasso}$ have the smallest meantime in combination with the highest rand index. The rand indices for $c_{Ridge}$ and $c_{Lasso}$ are $97.9\%$ and $98.6\%$ respectively. The lowest F1-score of $44.4\%$ is obtained by $c_{L1}$ which also has a high meantime. By using \texttt{PELT}, most predictions become more accurate with smaller meantime and higher F1-scores. Similar to the results for \texttt{WIN}, the results in Table~\ref{tab:res_PELT_simulated_datasets} show the highest F1-scores and rand index as well as
the smallest meantime error for $c_{Ridge}$ and $c_{Lasso}$. All cost functions except $c_{Normal}$ and $c_{AR}$ obtain F1-score of $100\%$, where the meantime vary between the algorithms. In Figure~\ref{fig:res_exp_data} (bottom) we see how $c_{L1}$, $c_{LinReg}$ and $c_{Ridge}$ give approximate indications of the true change points. Using the Bayesian approach, we get indications of the true change points and predicted change points after the exponential decay. Table~\ref{tab:res_BAYES_simulated_datasets} shows a meantime of $31.4$ time steps, F1-score $76.9\%$ and rand index $92.8\%$, which are all comparable to the predictions made using the optimisation approach.  

\subsubsection{Oscillation decay data}
\begin{figure}[h]
    \hspace{-0.7cm}
    \includegraphics[width = 1.05\linewidth]{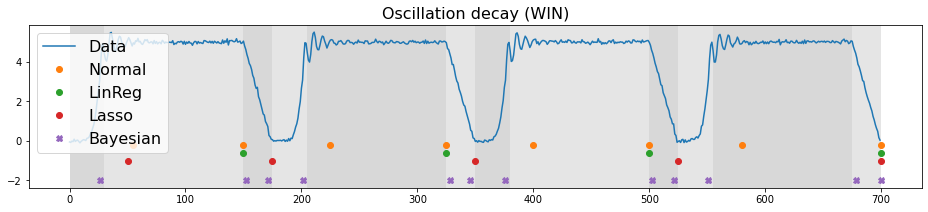}

    \hspace{-0.7cm}
    \includegraphics[width = 1.05\linewidth]{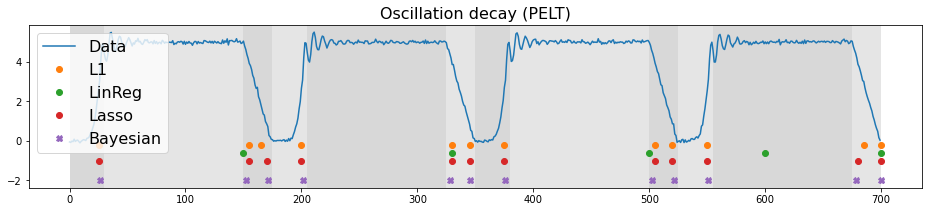}
    \caption{Results for the dataset with constant mean and changing variance, using various cost functions. Predicted change points are seen in coloured points, which are time points and do not have a value. The cost functions $c_{Normal}$, $c_{LinReg}$ and $c_{Lasso}$ can be seen in the upper figure, where the window-based search directions is used. The bottom image gives a comparison of $c_{L1}$, $c_{LinReg}$ and $c_{Lasso}$ when \texttt{PELT} is used. The predictions using \texttt{PELT} resemble the predictions made by the Bayesian approach, where three additional change points are predicted. }
    \label{fig:res_oscil_data}
\end{figure}
In this section we analyse data when oscillations can be found in process data instead of distinct changes in value, where the damping oscillations can obscure true change points. Figure~\ref{fig:res_oscil_data} shows a representation of a dataset with oscillations presented along with the predictions made using \texttt{WIN}, \texttt{PELT} and the Bayesian approach, respectively. The predictions presented in the bottom image agree on most change points. Table~\ref{tab:res_WIN_simulated_datasets} gives the metrics for the predictions, when we use the approximate search direction \texttt{WIN}, where number of predictions vary between the cost functions. The highest rand index is obtained using $c_{Lasso}$, $c_{Ridge}$ and $c_{AR}$ which all predict five change points. In exception to $c_{LinReg}$, $c_{Lasso}$ has the lowest meantime and a rand index of $89.1\%$. 
Using Table~\ref{tab:res_PELT_simulated_datasets} we observe that applying the optimal approach \texttt{PELT}, the F1-score and rand index are higher compared to the results for \texttt{WIN}.
The functions $c_{L2}$, $c_{L1}$, $c_{Normal}$, $c_{AR}$ and $c_{Lasso}$ predict twelve change points and have the highest rand indices of $97.0$, $97.1$, $95.6$, $96.8$ and $97.4$ percent respectively. The highest F1-score $81.8\%$ is obtained with $c_{Normal}$. Using the Bayesian approach, twelve change points are predicted, with a meantime of $3.5$ time steps, which is smaller than any of the predictions made by the optimisation approach. Similarly the rand index of $98.3\%$ is the highest of all predictions. 

\subsection{Real dataset}
\begin{table}[h]
    \centering
    \begin{tabular}{l c | c c | c c c c c}
        Cost & Penalty & K & AE & MT & Precision & Recall & F1 & RI \\
        function &&&& [seconds] & [\%] & [\%] & [\%] & [\%] \\
        \hline
        $c_{L2}$ & 3 & 12 & 5 & 67.9 & 72.7 & 50.0 & 59.3 & 96.4 \\
        $c_{L1}$ & 6 & 16 & 1 & 88.6 & 66.7 & 62.5 & 64.1 & 97.2 \\
        $c_{Normal}$ & 300 & 17 & 0 & 841.1 & 37.5 & 37.5 & 37.5 & 91.2 \\
        $c_{LinReg}$ & 6 & 16 & 1 & 194.0 & 60.0 & 56.3 & 58.1 & 92.9 \\
        $c_{AR}$ & 0.0015 & 15 & 2 & 144.8 & 50.0 & 43.8 & 46.7 & 96.5 \\
        $c_{Ridge}$ & 100 & 16 & 1 & 144.7 & 54.4 & 50.0 & 51.6 & 95.7 \\ 
        $c_{Lasso}$ & 100 & 17 & 0 & 135.9 & 50.0 & 50.0 & 50.0 & 95.8 \\ 
        \hline
    \end{tabular}
    \caption{Prediction results for PRONTO data using the optimisation approach, where the search method \texttt{WIN} is used. Each cost function indicated the best possible penalty level, along with respective obtained scores. }
    \label{tab:PRONTO_WIN}
\end{table}

\begin{table}[h]
    \centering
    \begin{tabular}{l c | c c | c c c c c}
        Cost & Penalty & K & AE & MT & Precision & Recall & F1 & RI \\
        function &&&& [seconds] & [\%] & [\%] & [\%] & [\%] \\
        \hline
        $c_{L2}$ & 150 & 16 & 1 & 191.9 & 66.7 & 62.5 & 64.5 & 96.4 \\
        $c_{L1}$ & 250 & 16 & 1 & 266.9 & 66.7 & 62.5 & 64.5 & 95.9 \\
        $c_{Normal}$ & 4500 & 23 & 6 & 294.8 & 54.5 & 75.0 & 63.2 & 96.1 \\
        $c_{LinReg}$ & 150 & 20 & 3 & 151.0 & 63.2 & 75.0 & 68.6 & 97.3 \\
        $c_{AR}$ & 0.02 & 22 & 5 & 343.4 & 38.1 & 50.0 & 43.2 & 91.6 \\ 
        $c_{Ridge}$ & 250 & 17 & 0 & 102.1 & 56.3 & 56.3 & 56.3 & 96.4 \\ 
        $c_{Lasso}$ & 250 & 17 & 0 & 99.9 & 56.3 & 56.3 & 56.3 & 96.4 \\ 
        
        \hline
    \end{tabular}
    \caption{Prediction results for PRONTO data using the optimisation approach, where the search method \texttt{PELT} is used. Each cost function indicated the best possible penalty level, along with respective obtained scores.}
    \label{tab:PRONTO_PELT}
\end{table}

\begin{figure}[h]
   \hspace{-0.7cm}
    \includegraphics[width = 1.05\linewidth]{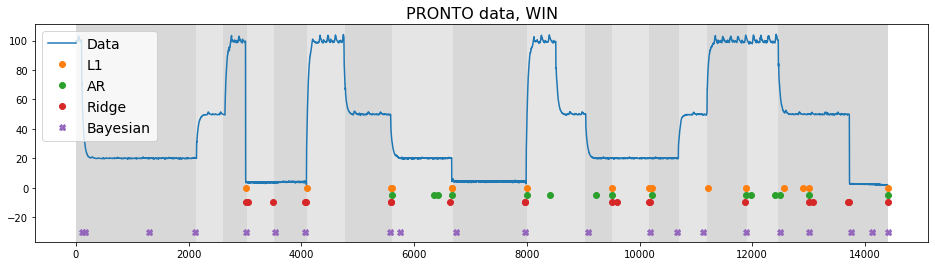}
    
    \hspace{-0.7cm}
    \includegraphics[width = 1.05\linewidth]{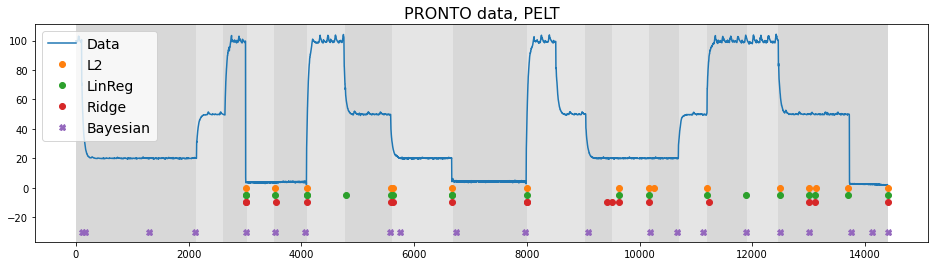}
    \caption{Results for the PRONTO dataset, using various cost functions. The signal is \texttt{Air~In1}, which is one out of four process variables presented in Figure \ref{fig:PRONTO_data}. Predicted change points are seen in coloured points, which are time points and do not have a value. The cost functions $c_{L1}$, $c_{LinReg}$ and $c_{Ridge}$ can be seen in both figures, where $c_{Ridge}$ agrees with the true change points using both \texttt{WIN} and \texttt{PELT}. The predictions resemble the predictions made by the Bayesian approach, where three additional change points are predicted. }
    \label{fig:res_PRONTO}
\end{figure}
We have tested different CPD algorithms for experimental datasets provided in \cite{PRONTOdataset}. The real world data presented in the PRONTO paper \cite{PRONTO_description} includes four separate process variable signals. The individual signals are seen in Figure~\ref{fig:PRONTO_data}, along with the change points present in all the process variables. The dataset is open source and can be retrieved via Zenodo \cite{PRONTOdataset}. Note, the signals are processed individually (uni-variate) and aggregated to get the combined prediction. Table~\ref{tab:PRONTO_WIN} shows the predictions obtained using \texttt{WIN}. To calculate the precision and recall, the margin $1\%$ of the number of samples are used, equivalent to $144 $ seconds in error is accepted as an accurate indication. The cost functions $c_{Normal}$ and $c_{Lasso}$ predict $17$ change points and have the lowest absolute error, where these algorithms also have the highest meantime. The smallest meantime is obtained by $c_{L2}$ and $c_{L1}$, which predict fewer points. The cost functions $c_{L2}$ and $c_{L1}$ have the highest F1-scores of $59.3\%$ and $64.1\%$ respectively and rand index of $96.4\%$ and $97.2\%$ respectively. The highest precision is obtained using $c_{L2}$ while the recall is lower at $50\%$. In Figure~\ref{fig:res_PRONTO} (top) we see a comparison of the predictions made by $c_{L1}$, $c_{AR}$ and $c_{Ridge}$. We see similarities in the predictions made by $c_{L1}$ and $c_{Ridge}$ and other predictions made by $c_{AR}$.
\\
\\
In the bottom image, we see predictions made by $c_{L2}$, $c_{AR}$ and $c_{Ridge}$ when \texttt{PELT} is used. The predictions' metric results are seen in Table~\ref{tab:PRONTO_PELT}. We can note the generally high rand index and F1-scores, with the exception of $c_{Normal}$. Most cost functions predict $17$ change points, and have a lower meantime compared to the values in Table~\ref{tab:PRONTO_WIN}. The smallest meantime of $95.8$ time steps is obtained by $c_{AR}$ and the largest meantime by $c_{Normal}$. The highest F1-score of $66.7\%$ and rand index of $95.4\%$ is obtained by $c_{Ridge}$. In Figure~\ref{fig:res_PRONTO} (bottom) we see that $c_{L2}$, $c_{LinReg}$ and $c_{Ridge}$ give similar predictions and correspond to many of the true change points. 
\\
\\
Figure~\ref{fig:res_PRONTO} also shows predictions obtained by the Bayesian method. The method predicts $21$ change points, which gives an absolute error $AE = 4$. The meantime of the predictions is $448.9$ seconds. The precision and recall are $65.0\%$ and $81.3\%$ respectively, which gives an F1-score of $72.2\%$. The rand index is $96.0\%$.
\\
\\
To give some insight to whether the regularisation parameter $\gamma$ in equations \eqref{eq:c_ridge}~and~\eqref{eq:c_lasso} affects the predictions, different parameter values can be chosen while all other parameters are unchanged. Tables \ref{tab:PRONTO_reg_ridge} and \ref{tab:PRONTO_reg_lasso} present the predictions made using the two cost functions $c_{Ridge}$ and $c_{Lasso}$ respectively. Both cost functions are applied using search direction \texttt{WIN} and with a penalty term $pen = 100$. Both tables show how the regularisation parameter can influence the predictions, where a higher regularisation parameter gives fewer predictions. We can notice how the meantime is reduced when the parameter value increases. 
\begin{table}[H]
    \centering
    \begin{tabular}{l | c c | c c c c c}
        Reg. constant & K & AE & MT & Precision & Recall & F1 & RI \\
        ($\gamma$) &&& [seconds] & [\%] & [\%] & [\%] & [\%] \\
        \hline
        0.1 & 16 & 1 & 144.7 & 53.3 & 50.0 & 51.6 & 95.7 \\
        1 & 16 & 1 & 144.7 & 53.3 & 50.0 & 51.6 & 95.7 \\
        10 & 16 & 1 & 144.7 & 53.3 & 50.0 & 51.6 & 95.7 \\
        100 & 16 & 1 & 144.7 & 53.3 & 50.0 & 51.6 & 95.7 \\
        1000 & 17 & 0 & 136.6 & 50.0 & 50.0 & 50.0 & 95.8 \\
        10000 & 12 & 2 & 178.9 & 54.5 & 37.5 & 44.4 & 93.1 \\
        \hline
    \end{tabular}
    \caption{Prediction results for PRONTO data using the optimisation approach and different regularisation constant $\gamma$ in $c_{Ridge}$, see equation \eqref{eq:c_ridge}. The search method \texttt{WIN} is used and a penalty $pen = 100$. We can note that different predictions are made depending on which regularisation constant is used. The predictions in Tables \ref{tab:PRONTO_WIN} and \ref{tab:PRONTO_PELT} use $\gamma = 1$, and are included in the table to be used as comparison.}
    \label{tab:PRONTO_reg_ridge}
\end{table}

\begin{table}[H]
    \centering
    \begin{tabular}{l | c c | c c c c c}
        Reg. constant & K & AE & MT & Precision & Recall & F1 & RI \\
        ($\gamma$) &&& [seconds] & [\%] & [\%] & [\%] & [\%] \\
        \hline
        0.1 & 17 & 0 & 284.3 & 68.8 & 68.8 & 68.8 & 94.5 \\
        1 & 17 & 0 & 135.9 & 50.0 & 50.0 & 50.0 & 95.8 \\
        10 & 17 & 0 & 135.9 & 50.0 & 50.0 & 50.0 & 95.8 \\
        100 & 7 & 10 & 104.3 & 83.3 & 31.3 & 45.5 & 87.7 \\
        1000 & 8 & 9 & 97.3 & 71.4 & 31.1 & 43.5 & 94.5 \\
        10000 & 8 & 9 & 97.3 & 71.4 & 31.1 & 43.5 & 94.5 \\
        \hline
    \end{tabular}
    \caption{Prediction results for PRONTO data using the optimisation approach and different regularisation constant $\gamma$ in $c_{Lasso}$, see equation \eqref{eq:c_lasso}. The search method \texttt{WIN} is used and a penalty $pen = 100$. We can note that different predictions are made depending on which regularisation constant is used. The predictions in Tables \ref{tab:PRONTO_WIN} and \ref{tab:PRONTO_PELT} use $\gamma = 1$, and are included in the table to be used as comparison.}
    \label{tab:PRONTO_reg_lasso}
\end{table}


\section{Discussion}
\label{sec:discussion}
In this section we analyse how and why different  methods for CPD are chosen as well as discuss results presented in sections \ref{sec:methods} and \ref{sec:results}, respectively. First, we discuss the testing procedure and motivate some parameters used in this work. Then we analyse obtained results for the two approaches presented in section \ref{sec:results}. Finally, the user interaction of each approach is discussed, followed by suggestions for future work. 

\subsection{Testing procedure}
All predictions are made on uni-variate signals. This is because the used implementation for the Bayesian approach is not able to make predictions on multi-variate datasets. The optimisation approach is able to make simultaneous predictions based on multiple correlated signals, which in some cases give other predictions. 
Using our numerical investigations we can conclude that
this happens because more information is incorporated into the prediction algorithms. Therefore, some change points are detected in the multi-dimensional case and not in the uni-variate case. 
To make the approaches more comparable, predictions are made only using the uni-variate signals and aggregating the results for the signals.
This is mainly relevant for the predictions made on the PRONTO dataset, where we have four process variables. The predictions are selected using different methods for the two approaches. For the optimisation approach, the penalty term is tuned, while the peak thresholds can be tweaked for the Bayesian approach.
\\
\\
Another aspect to discuss is the choice of penalty in the optimisation approach. 
A benefit with the approach is that the desired number of change points is not known a priori; we also do not use the number of predictions as a parameter in the algorithm. The algorithm is applied multiple times with different penalty values, and the best prediction is selected based on metric values. This can be computationally heavy if many iterations are necessary, but when the appropriate penalty level is found, it can be used for other subsets of the same data. 
We want to minimise AE and meantime, while we want to maximise the F1-score and rand index. Some extra emphasis has also been given to precision over recall, where we value getting an accurate prediction over covering all true change points. 
This is on the basis that some change points might not be identifiable in the dataset but can be of specific interest, according to Definition \ref{def:CP_domain}. Using the F1-score as a function of penalty we will generally find a maximum value. This is since on either side of the maximum, either the recall or precision are lowered due to the number of predictions diverging from the true number of change points. The maximum value for the F1-score has generally been selected as the best prediction.
It should be noted that the results presented in the tables above are not necessarily the only good predictions, and another prediction may be selected if different metrics are used. 
\\
\\
For the Bayesian approach, predictions are selected based on sampling from the posterior distribution. The selection is done according to the Maximum A Posteriori (MAP) approach such that the points with the maximum posterior probability are selected as the detected change points. In the implementation, the peaks of the posterior distribution are identified and returned as change points. To be identified as a peak, the point has to be a local maximum, with value above a threshold and not too close to another change point. These values can be alternated by the user to change some of the predictions. The confidence level used as a threshold for a peak is set to $50\%$ and the a proximity threshold of $10$ time steps for each dataset. The confidence level is based on that we want the probability of being a change point being majority compared to the probability of not being a change point. Similar to the optimisation approach, the number of selected points is not known a priori, and the approach is not dependent on a specified number of change points. 
This underlines that changing the peak parameters will not necessarily give new predictions or give the exact number of change point predictions. 
It should be noted that other sampling methods can be applied to the same posterior distribution to get other predictions. 
%

\subsection{Test metrics}
The metrics are chosen to give a detailed overview of the performance from different perspectives. The meantime indicates the best case scenario, where it measures the average time between predictions and the closest true change point. This, however, does not take into account if the distance for a prediction is measured to the corresponding true change point, or merely the closest one. The precision and recall indicate the accuracy in the predictions and coverage of true change points respectively. 
The metrics can be combined to an average metric, in this case the F1-score. The F1-score is the harmonic mean between the two metrics, but this can be chosen as some other metric. 
In addition to the harmonic mean, if precision is more important than recall, the combined metric could be a weighted average instead. 
In this case, the F1-score is chosen as it weighs the metrics equally since we do not know the respective importance levels of the two metrics. 
However, 
sometimes precision is given more credit than recall when selecting the best model. 
In contrast to only evaluating the predicted change points, the rand index compared the obtained segmentation, and therefore gives a higher level comparison of the overall prediction. 
This metric indicates both the amount of agreements, and also the amount of disagreements.
In addition to the used metrics, other metrics such as network based distances can be used depending on desired information. 
Examples are the Hausdorff metric which measures how far the predicted and actual change points are from each other. Another network based metric is the Hamming distance, which measures the minimum number of substitutions required for the predicted and actual set of change points to be identical. 
\\
\\
The choice of metrics is of extra interest when the choice of predictions are based on them, and some of the metrics can influence each other. 
Naturally, the number of predictions indicate how many changes were identified, but does not give information on how accurate the predictions are. 
The meantime gives the average distance between the predictions and the closest true change point, which give an indication of the accuracy. 
On the other hand, the meantime will be affected by the number of predictions, where we can expect a lower meantime if we predict fewer change points. 
An example of this is seen in Table~\ref{tab:PRONTO_WIN}, where $c_{L2}$ has a meantime of $67.9$ seconds and an F1-score of $59.3\%$, while $c_{Lasso}$ has a meantime of $284.3$ seconds and F1-score $68.8\%$. 
We can note that the number of predictions are lower for $c_{L2}$ and the precision is higher compared to $c_{Lasso}$. 
Similarly, we can also note that the rand index is related to the meantime. With the lower meantime, $c_{L2}$ has a higher rand index than $c_{Lasso}$, which has a higher meantime. 

\subsection{Results}
In this section we discuss the obtained results from the various approaches. 
All predictions are made on uni-variate datasets, where six simulated datasets are used as well as a real world dataset consisting of four process variables. 
To make fair predictions between the process variables, these four process variables are normalised before prediction. 
Normalisation could also be applied to the simulated datasets, bu our computations show that it does not affects the predicted change points, merely the value of the penalty term. 
In our computational results, all actual change points are not identified, especially for the real world dataset. 
This is since some change points are domain specific, according to Definition \ref{def:CP_domain}, and do not have a distinct change in the data linked to the change point. 
We cannot expect the algorithms to detect these change points as they are designed to identify points defined in Definition \ref{def:CP}. 
In the following sections, the simulated and real world datasets are discussed separately. 

\subsubsection{Simulated datasets}
\begin{figure}[H]
    \hspace{-0.7cm}
    \includegraphics[width = 1.05\linewidth]{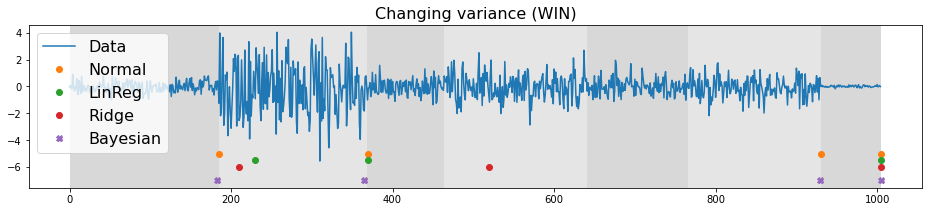}

    \hspace{-0.7cm}
    \includegraphics[width = 1.05\linewidth]{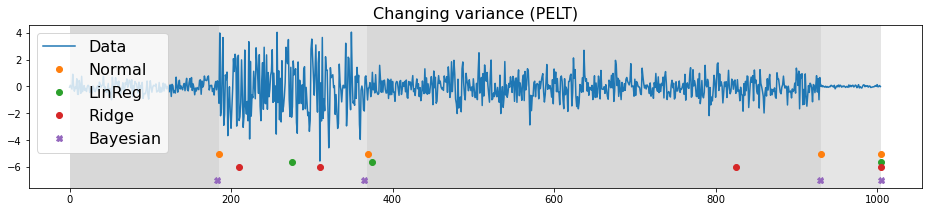}
    \caption{Illustration of results  for  the  dataset  with  constant  mean  and  changing  variance with fewer predictions than indicated in Figure~\ref{fig:res_changing_var}. In this case we see more similarity between the predictions where only four change points are predicted.}
    \label{fig:res_changing_var_special}
\end{figure}

The change points in the piecewise constant dataset are identified by most cost functions in the optimisation approach and by the Bayesian approach. 
These types of features are common in processes, and we see piecewise constant signals in the PRONTO dataset (see Figure~\ref{fig:PRONTO_data}).
In Tables \ref{tab:res_WIN_simulated_datasets} and \ref{tab:res_PELT_simulated_datasets} we see that $c_{L2}$, $c_{Normal}$, $c_{Ridge}$ and $c_{Lasso}$ manage to detect six out of seven change points with high accuracy. 
The maximum likelihood based cost functions $c_{L2}$ and $c_{Normal}$ identify the change in mean as a large cost, according to equation \eqref{eq:c_l2} and \eqref{eq:c_normal_univar} respectively. 
Using $c_{L1}$ we also get good predictions, but in this case the meantime is higher, and accuracy lower. 
This happens since it takes deviation into account and not variance.
The model fitting approaches, such as $c_{LinReg}$ and $c_{AR}$, generally have a lower recall than other approaches, indicating that all true change points are not predicted accurately. 
Fitting an autoregressive model gives a lower meantime than using ordinary linear regression, suggesting that incorporating previous samples improves the predictions. 
Using Tikhonov regularisation, as we have done in $c_{Ridge}$ and $c_{Lasso}$, also gives good predictions, where especially $c_{Lasso}$ gives a lower meantime. 
Since these cost functions are extensions of ordinary linear regression, see equations \eqref{eq:c_ridge} and \eqref{eq:c_lasso}, a comparison can be made between $c_{LinReg}$ and the regularised cost functions.
Both regularised cost functions improve the results compared to solely linear regression. 
Improvements are mainly seen in meantime, recall and rand index. 
The improvements come from the regularised cost functions not trying to over-fit to the variability in the data, which should give a penalty in the regularisation term, whereas this additional penalty is not present in the ordinary linear regression. 
However, looking at Figure~\ref{fig:res_piecewise_constant}, we see that the regularisation approach gives some accurate predictions, while some change points are predicted to give segments where a regularised linear model can be fitted over multiple true segments, note especially interval $t_i \in [800, 1200]$ for \texttt{PELT}. 
The Bayesian approach manages to predict all seven points quite accurately. 
The meantime is below $10$ time steps, which is somewhat higher than the median meantime for the optimisation approach, while the precision, recall and rand index are equal to the best prediction methods of the optimisation approach. 
This indicates that the Bayesian approach with a flat prior and Gaussian likelihood function is suitable, and gives good predictions, for the piecewise constant dataset. 
\\
\\
If the piecewise constant segments are changed to piecewise linear segments, some approaches struggle to predict the change points accurately. 
In this setting, it can be presumed that model based cost functions in the optimisation approach will give better predictions, compared to maximum likelihood cost functions. 
We can also note that the range of the data in this set, $y_{t_i} \in [-1, 1]$, is the smallest of the simulated datasets. 
In Table~\ref{tab:res_WIN_simulated_datasets} we see how the cost functions $c_{Normal}$, $c_{AR}$ and $c_{Ridge}$ give the best predictions when \texttt{WIN} is used. 
The latter two are model based functions, which are assumed to make better predictions, while the normal cost function is a maximum likelihood function which incorporates more information than the norm based functions $c_{L2}$ and $c_{L1}$. Using the regularised cost function, the meantime is zero and F1-score $100\%$, suggesting that it predicts all the change points perfectly. 
When the exact search direction is used, the results in Table~\ref{tab:res_PELT_simulated_datasets} are improved compared to the approximate search direction. 
The number of predictions has increased compared to \texttt{WIN}, which can influence the meantime and precision negatively, while increasing the recall. 
In the case of $c_{AR}$, we see the same number of predictions in the two tables, and the meantime is somewhat lower when \texttt{PELT} is used. 
On the other hand, the precision and recall are lower. 
This suggests that in general the predictions are more accurate, while one prediction is outside the allowed margin for precision and recall, causing a decrease in the meantime. The Bayesian approach predicts only two points, where one is the artificial change point at the end of the data.
Illustrations are presented in Appendix A in \cite{masterThesis}.
\\
\\
If only the variance changes, then it is natural to expect approaches which incorporate the variance to give the best predictions. The Bayesian approach includes information concerning change in variance in the Gaussian likelihood function, where a change in variance can mean a change in distribution and therefore, a change point. The maximum likelihood based cost functions $c_{L2}$ and $c_{L1}$ do not incorporate the variance directly, while the extended version $c_{Normal}$ does, see equation \eqref{eq:c_normal_univar}. In Table~\ref{tab:res_WIN_simulated_datasets} we see how this cost function gives the most accurate predictions, with $AE = 0$ and a meantime of $22$ time steps, which is the lowest of all cost functions. We can also note how $c_{Lasso}$ does not identify any change points (even with zero penalty) in this setting, except the artificial one at the end, which suggests that it fits one model for the entirety of the data. In Figure~\ref{fig:changing_var_dataset} (top) we see how $c_{Normal}$ predicts four of the change points accurately, while the remainder are not accurate predictions. This can arguably be due to the features in the data. Some change points do not have as distinct change points, as defined in Definition \ref{def:CP}. The lack of distinctness in the change points causes the change point detection to be difficult, as the change is not significant. In Figure~\ref{fig:changing_var_dataset} we also see that this correlates with the predictions made by the Bayesian approach, where four change points are identified. To extend the investigation, the algorithms can predict fewer change points and see if these agree more. In Figure~\ref{fig:res_changing_var_special} we see the predictions made by a few algorithms and the Bayesian approach. The predictions made using the Bayesian approach and $c_{Normal}$ give similar predictions, while other cost functions give different predictions. This confirms that both the Bayesian approach and the $c_{Normal}$ cost function incorporate variance changes in an effective way, while other cost functions do not. 
\\
\noindent
Autoregressive data, where values are linked to previous values, is common in many areas and therefor the dataset similar to the one seen in Figure~\ref{fig:AR_dataset} is relevant to investigate further. In this case, we only have one approach which incorporates the autoregressive principle, namely $c_{AR}$, which is presumed to give the best predictions. This is confirmed in Tables \ref{tab:res_WIN_simulated_datasets} and \ref{tab:res_PELT_simulated_datasets} where $c_{AR}$ gives the best predictions, especially when \texttt{WIN} is used. Maximum likelihood based cost functions cluster the predictions around the areas with the highest variability, while regularised functions fits fewer models and predicts fewer change points (with \texttt{WIN}). The linear regression model $c_{LinReg}$ predicts the correct number of change points, but do not predict the true change points, see Figure~\ref{fig:res_AR_data}. The autoregressive model predicts all change points correctly, as it incorporates the previous sample in the model fitting, which is not the case for the other cost functions. This model fitting works especially well when \texttt{WIN} is used, due to using a window view and not taking the entirety into account as in \texttt{PELT}, which reduces the risk of over-fitting. The Bayesian approach predicts similar change points as the maximum likelihood, where distinct changes in the data (especially in segments with high variability) are perceived as change points. This suggests that the Bayesian approach does not take information from previous samples into account. 
\\
\\
The exponential decay dataset has similarities to real-world processes, where the exponential decay is a common phenomenon, this particular feature is of special interest. 
When \texttt{WIN} is used, the cost functions $c_{Normal}$, $c_{Ridge}$ and $c_{Lasso}$ manages to to get F1-scores of $100\%$ and RI above $95\%$, where $c_{Normal}$ has the highest meantime, see Table~\ref{tab:res_WIN_simulated_datasets}. 
When \texttt{PELT} is used, in addition to the three mentioned cost functions, $c_{L2}$ and $c_{L1}$ manages to get F1-score of $100\%$, with a higher meantime than $c_{Ridge}$ and $c_{Lasso}$ but lower than $c_{Normal}$. 
This suggests that using a regularised model gives the best predictions with the lowest meantime of below ten time steps for both search directions. 
Even though there is variance present and different trends in the data, the cost functions in equations \eqref{eq:c_ridge} and \eqref{eq:c_lasso} fits models to the various segments without over-fitting to the variability. 
The drastic improvements in the predictions by $c_{L2}$ and $c_{L1}$ when \texttt{PELT} is used lies in the limited number of predictions made using \texttt{WIN}. 
When the approximate search direction is used without additional penalty ($pen = 0$), the algorithms still only predicts $4-5$ change points. 
Using the exact search direction, more partitions are investigated, and more points can be detected, see Table~\ref{tab:res_PELT_simulated_datasets}. 
Using Table~\ref{tab:res_BAYES_simulated_datasets}, we observe that the Bayesian approach manages to predict nine change points with $100\%$ precision, and a rand index of $92.8\%$.
\\
\\
Similar to the exponential decay data, the oscillating phenomenon seen in Figure~\ref{fig:oscil_data} are common in processes when a new level stabilises after an alternation is imposed. 
The number of predictions when \texttt{WIN} is used is significantly lower then when \texttt{PELT} is used, see Tables \ref{tab:res_WIN_simulated_datasets} and \ref{tab:res_PELT_simulated_datasets}. 
In this case, it is more reasonable to use \texttt{PELT} as the dataset is not too large and the time complexity of the calculations is low. 
Then, the number of predicted change points corresponds better to the number of true change points. 
In Figure~\ref{fig:res_oscil_data} (bottom) we see the predictions of $c_{L1}$, $c_{LinReg}$, $c_{Lasso}$ and the Bayesian approach, where most approaches give similar predictions. 
Visually, the best predictions are made by the Bayesian approach, which is also confirmed in the metrics presented in Table~\ref{tab:res_BAYES_simulated_datasets} where the meantime is lower than any predictions presented in Table~\ref{tab:res_PELT_simulated_datasets}. 
Illustrations are presented in Appendix A in \cite{masterThesis}.
This shows that the Bayesian method is appropriate to use when oscillations as well as various linear 
trends are present in the signal. 
This brings us further to the study of the real world dataset, where these phenomenons are present. 

\subsubsection{Real dataset}
The real world dataset provided by PRONTO incorporates some of the features seen in the simulated datasets. 
In addition to this, the dataset is larger than the simulated datasets, with approximately $14000$ samples. 
In Figure~\ref{fig:PRONTO_data} we see the four process variables \texttt{Air In1}, \texttt{Air In2}, \texttt{Water In1} and \texttt{Water In2} studied in this work. 
In all variables we observe piecewise constant segments as well as some exponential trends between segments. We also see variability, especially in \texttt{Air In2}. 
In addition, we can presume that there exists some variance present in the entirety of the signals, and we can note that the predictions are made on normalised signals. 
This indicates partial similarities with the piecewise constant, exponential decay and changing variance datasets. 
Investigating the performance of cost functions in the mentioned simulated dataset, we find that generally $c_{Normal}$ performs quite well in all cases, and the regularised approaches $c_{Ridge}$ and $c_{Lasso}$ perform well in the piecewise constant and exponential decay dataset. 
An assumption would be that these cost functions will give good predictions, based on the trials on the simulated datasets.
\\
\begin{figure}[h]
    \hspace{-0.7cm}
    \includegraphics[width = 1.05\linewidth]{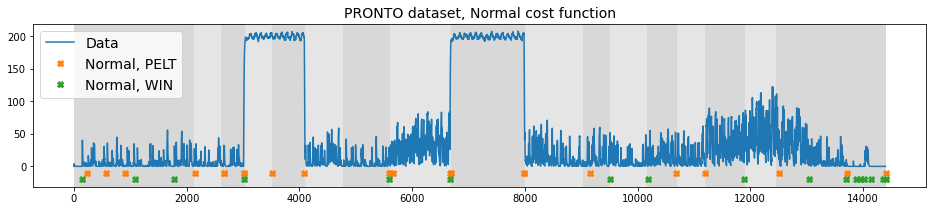}
    \caption{Results for the PRONTO dataset, using $c_{Normal}$ with \texttt{PELT} and \texttt{WIN}. The signal is \texttt{Air~In2}, which is one out of four process variables. We can notice how not all the change points are identified. }
    \label{fig:res_PRONTO_normal}
\end{figure}
\\
Tables \ref{tab:PRONTO_WIN} and \ref{tab:PRONTO_PELT} show that $c_{Normal}$
does not give the best change point predictions.
To explore this further, Figure~\ref{fig:res_PRONTO_normal} shows the predictions made by $c_{Normal}$ using both \texttt{PELT} and \texttt{WIN}. 
Note that the signal in the background is the \texttt{Air In2} since this variable shows the changes in variance explicitly. 
In the figure, we see how the approximate solution indicates change points where the variance changes, but which is not necessarily a change point, see interval around $t_i = 14000$ where multiple detections are found due to irregular variance. 
When \texttt{PELT} is used, the algorithm predicts $23$ change points, which gives an annotation error of six, which is the highest of all cost functions. 
In Figure~\ref{fig:res_PRONTO_normal} we see how some of the predictions are clustered in the first $1000$ time steps, where no change point is present. 
The conclusion from this illustration is that $c_{Normal}$ might emphasise on features which are not necessarily linked to change points, in some sense over-fitting to the variance. 
This is also indicated in the predictions made by the simpler cost function $c_{L2}$, which gives better predictions than $c_{Normal}$. 
Note how $c_{Normal}$ is an extension of $c_{L2}$ which incorporates the variance, see equation \eqref{eq:c_normal_univar}. 
The function $c_{L2}$ gives good predictions on the piecewise constant dataset, which are identical to the ones made by $c_{Normal}$. 
This suggests that the variance changes are not the main features linked to change points in the PRONTO dataset, and where $c_{Normal}$ over-fits to the variance. 
To oppose this over-fitting, we can use the regularised cost functions. 
\\
\\
Using the regularised cost function we get the best predictions, along with using $c_{L2}$. 
The metrics are seen in Tables \ref{tab:PRONTO_WIN} and \ref{tab:PRONTO_PELT}, where F1-score is increased when \texttt{PELT} is used instead of \texttt{WIN}. 
The cost function $c_{L2}$ has better precision and recall than the regularised methods, but the annotation error and meantime are reduced when the regularised cost function is used with \texttt{PELT}. 
Using the regularised cost functions give better predictions than $c_{Normal}$. 
This suggests that over-fitting to the variance should be avoided in this case, where $c_{Normal}$ accounts for changes which are not necessarily linked to change points. 
For this dataset, using the regularised cost function in combination with \texttt{PELT} gives the best predictions. 
\\
\\
The Bayesian approach has not shown tendencies to over-fit predictions in the simulated dataset, rather under-fit. 
In Figure~\ref{fig:res_PRONTO} we see how the $21$ predicted change points cover most of the change points, and also give some additional indications. 
We can note that the recall is the highest out of all predictions, but also the meantime. 
The significantly higher meantime is due to the multiple predictions made in the beginning and end of the dataset, which are places far away from the closest change point. 
These indications are due to the significant drop in the signal \texttt{Air in1} at around $t_i=100$ and the significant change in variance in \texttt{Air In2} at around $t_i = 14000$. 
Accordingly to Definition \ref{def:CP}, these points correspond to significant changes in data, and we are expecting to get a change point  here in the case when they will be sensitive to changes in mean or variance.
In this case the Bayesian approach incorporates both. 
If the predictions in the beginning and end of time intervals would have been correct, the predictions of the Bayesian method would yield a good prediction of most of the true change points present in the PRONTO dataset. 
\\
\\
In Tables \ref{tab:PRONTO_reg_ridge}~and~\ref{tab:PRONTO_reg_lasso} we see how different regularisation parameters can affect the prediction. 
The tables show that selecting an appropriate regularisation constant is of importance. In the results presented in Table~\ref{tab:res_WIN_simulated_datasets}-\ref{tab:PRONTO_PELT} the constant $\gamma=1$ has been used, where Tables \ref{tab:PRONTO_reg_ridge}~and~\ref{tab:PRONTO_reg_lasso} present other possible values for $\gamma$. 
It is evident that the investigation is not exact, and indicates that different parameter choices for $\gamma$ give different predictions. 
In contrast to simply selecting different values, other methods to estimate an appropriate value can be used, such as iterative approximation. 
These extended algorithms are not covered in this work, but could be of interest for future work. 

\subsection{User interaction}
To make the explored unsupervised CPD approaches user friendly and more effective, we want to explore the possibilities of incorporating user feedback or prior knowledge. 
This can either be done via changes in settings of the unsupervised approaches or after predictions have been made. 
The two approaches have different modifications which can be made to alter the predictions. 
In the upcoming sections these possibilities are discussed with respect to the optimisation approach and the Bayesian approach individually. 

\subsubsection{Optimisation approach}
After the discussion of the results in previous sections, a natural way to affect the predicted change points is to select an appropriate cost function. If the data contains certain features, as explored in the simulated datasets, we can chose an appropriate cost function according to the underlying assumptions and function formulation. Alternatively, we know which type of change point we want to detect, and therefor can choose a cost function which is sensitive to specific changes in the data.
\\
\\
A parameter which can be changed to affect the number of predictions is the penalty. This constant can be estimated by performing multiple estimations and selecting the penalty level which gives the best metrics. If this was not sought after, the constant can be either increased to reduce the number of predictions or decreased to possibly\footnote{If the penalty is zero, the penalty term cannot be reduced further and therefor not generate more predictions.} increase the number of predictions. A user friendly implementation would first estimate the best penalty level and perform CPD. This penalty can then be altered by the user, if the results are not satisfactory, to give another prediction. 
\\
\\
If the unsupervised algorithm gives a generally good prediction of the change points, but there are domain specific change points present (see Definition \ref{def:CP_domain}), the predicted points can be altered. A change point could either be added or removed given the prediction made by the unsupervised algorithm. It should be noted that this will be a tedious task for a large number of domain specific change points, and additional solutions will be required to automate such tasks. 

\subsubsection{Bayesian approach}
For the Bayesian approach we can make appropriate choices for prior and likelihood function, similar to the choice of cost function in the optimisation approach. The prior contains information about the distance between change points while the likelihood function contains information about the segments. By choosing these according to the available data, the accuracy of the predictions can be increased. An example is to use a discrete geometric distribution as a prior, instead of a flat prior, if there is knowledge of the change points being equally spaced. 
\\
\\
A parameter change can be done to the peak selection criterion. As default, the peak selection can have a setting similar to the one used in this work, to give a first prediction of the change points. If the results are not satisfactory, the user could try to change the confidence level or distance parameters and see if this gives other predictions. If further alternations are necessary, a predictions can either be added or removed using the same principle as the optimisation approach. 
\\
\\
In addition to the presented user interaction, the use of probabilities can incorporate additional probability distributions. When the posterior distribution is calculated, but an area which should contain a change point does not have a high probability, the posterior distribution can be further joint with a domain specified distribution. This means that if a user knows approximately where a change point is present, a probability distribution with mean and variance according to the users' expectation can be found and joined with the calculated posterior distribution. This means that the posterior distribution of change points can incorporate the information from users without being recalculated.

\subsection{Future work}
As the application area of change point detection is wide and options of how to predict these points are many, this work does not cover all relevant topics. 
Future work can be pursued in areas such as defining new cost functions linked to specific feature identification using the optimisation approach, or other prior and likelihood function for the Bayesian approach. 
By making strategic choices, and trying to incorporate information about the data or sought after change points, the accuracy in the predictions can be increased. 
The application area can also be widened by studying the online version of change point detection, whereas this work has focused on the offline version. 
Potentially, similar comparisons can be made for the online CPD where aspects of this work can be implemented. 
\\
\\
This work has focused on simple implementations of Tikhonov regularisation, but there are numerous other options. This work has studied a selected constant as the regularisation parameter $\gamma$ in equations \eqref{eq:c_ridge}~and~\eqref{eq:c_lasso}, but there is potential in iterative updates of the parameter \cite{tikhonov_numerical_methods_for_ill_posed_problems, iter_regularizartion_for_ill_posed_problems}. Other regularisation methods which can be studied include the balancing principle \cite{num_analysis_LS_and_perceptron_for_classification}, L-curve or S-curve methods
\cite{niinimaki} and randomised SVD methods \cite{Ito}. Total variation regularisation can also be applied, which has been applied for denoising data \cite{total_variation}.
\\
\\
In addition to the theoretical aspects, the usability in industry settings can be studied further. This work has focused on comparing the usability in industries, where further investigations can be pursued in terms of change point usage. A relevant question is: how can the production in industrial processes benefit from the predicted change points in the process data? Ideas could be to use the change points to identify production phases or for syncing phases in chemical processes in order to identify anomalies.

\section{Summary}
\label{sec:summary}
This work has studied two unsupervised  algorithms
 based on two different approaches
for detecting change points.
 One approach is formulated as an optimisation problem while the other is based on Bayesian statistics. 
The study has found that the algorithms are affected by the features in the underlying data, where both simulated and real world data have been examined. 
In the optimisation approach, the choice of cost function can affect the predictions made by the algorithm, where these functions are strongly linked to the type of features being detected as change points. 
A new type of cost function has been introduced, which uses Tikhonov regularisation.
The Bayesian approach uses prior knowledge on the distance between change points and a likelihood function with information about the segments in order to predict the probability of a time point being a change point.  
To this day, the two mentioned approaches have not been compared, nor has the regularisation been used for change point detection previously. 
In addition to this, the possibility to incorporate user feedback has been explored, where both approaches are able to incorporate user knowledge and corrections post prediction.
This sheds light on the potential usability outside of academia.
The importance of change point detection becomes more and more important as datasets increase in size, where these unsupervised detection algorithms can help users process the data in order to draw conclusions.

\newpage
\bibliography{bib.bib}
\bibliographystyle{unsrt}

\end{document}